\documentclass{article}
\usepackage{amssymb}
\usepackage{amsmath}
\usepackage{amscd}

\title{On congruences between Drinfeld modular forms}
\author{Arash Rastegar}

\begin{document}

\maketitle

\begin{abstract}
Let ${\mathbf F}_q$ denote a finite field of characteristic $p$
and let $n$ be an effective divisor on the affine line over
${\mathbf F}_q$ and let $v$ be a point on the affine line outside
$n$. In this paper, we get congruences between ${\mathbb
Q}_l$-valued weight two $v$-old Drinfeld modular forms and
$v$-new Drinfeld modular forms of level $vn$. In order to do
this, we shall first construct a cokernel torsion-free injection
from a full lattice in the space of $v$-old Drinfeld modular
forms of level $vn$ into a full lattice in the space of all
Drinfeld modular forms of level $vn$. To get this injection we
use ideas introduced by Gekeler and Reversat on uniformization of
jacobians of Drinfeld moduli curves.
\end{abstract}

\section*{Introduction}

In the function field case, the information given by Hecke action
on Jacobians of compactified Drinfeld moduli curves is richer than
the Hecke action on the space of automorphic forms. As a
consequence, multiplicity one fails to hold mod-$p$ and Hecke
action will no longer be semi-simple. This points to the analytic
nature of these building blocks of the theory of elliptic modular
forms, which are absent in the function field case. In view of
this disorder, it is crucial to understand the relation between
geometric and analytic theories.

The results of Gekeler and Reversat confirm the importance of
understanding congruences between Drinfeld modular forms. They
prove that double cuspidal Drinfeld modular forms of weight two
and height one with ${\mathbb F}_p$ residues are reduction modulo
$p$ of automorphic forms. And only double cuspidal harmonic
cochains with ${\mathbb F}_p$ coefficients can be lifted to
harmonic cochains with $\mathbb Z$ coefficients [Ge-Re]. Finding
congruences between Drinfeld modular forms is also an important
step towards formulation and development of Serre's conjectures
on weight and level of modular forms. This is the main
application we emphasize on. We use uniformizations for the
Drinfeld moduli spaces introduced by Gekeler and Reversat to
obtain congruences between Drinfeld modular forms of rank $2$.
The computaions of the congruence ideal are similar to the number
field case due to Ribet [Ri]. We also calculate the congruence
module for more complicated congruences in the final section.
These congruences can be used to obtain towers of congruences
between Hecke algebras.

\section{Background on Drinfeld modular forms}\label{sec1}

In the number field case, the notion of modular form is based on
the moduli space of principally polarized abelian varieties. In a
series of papers, Drinfeld introduced the notion of Drinfeld
module as a function field analogue to an abelian variety and
defined the notion of Drinfeld modular form using the moduli
space of Drinfeld modules. He succeeded to prove a special case of
Langlands conjectures in the function field context (see [Dr1],
[Dr2]). Our main references in this chapter are [Dr3], [Ge],
[Ge-Re], [Go] and [Go-Ha-Ro].

\subsection{Drinfeld moduli spaces}\label{ssec1.1}

Let $X$ be a smooth projective absolutely irreducible curve of
genus $g$ over ${\mathbf F}_q$. The field of rational functions
$K$ of the curve $X$ is an extension of ${\mathbf F}_q$ of
transcendence degree one. We fix a place $\infty$ of $K$ with
associated normalized absolute value. We use the same notation
for the places of the field $K$ and ideals of its ring of
integers $A$ and effective divisors on $X$. Every place $v$ which
is not equal to $\infty$ is called a finite place. Let $q_v$
denote the order of the residue field of the ring of integers
$O_v$ of the completion $K_v$. The basic example will be the
function field ${\mathbf F}_q (t)$. The completion of the
algebraic closure of $K_{\infty}$, will respect to the unique
extension of absolute value is denoted by $C$.

Let $L$ be a field extension of ${\mathbf F}_q$ with an
$A$-algebra structure $\gamma : A \to L$ and $\tau$ denote the
endomorphism $x \mapsto x^q$ of the additive group scheme
underlying $L$. The endomorphism ring of the additive group
scheme underlying $L$ is the twisted polynomial ring $L \{ \tau_p
\}$ where $\tau_p : x \mapsto x^p$ satisfies the commutation rule
$\tau_p \circ x = x^p \circ \tau_p$ for all $x \in L$. Let $S$ be
a scheme over ${\rm Spec} \, A$, an $L$-Drinfeld module $({\mathbb
L} , \Phi)$ of rank $r \in {\mathbb N}$ over $S$ consists of a
line bundle ${\mathbb L}$ over $S$ and a ring homomorphism $\Phi
: A \to {\rm End}_S ({\mathbb L} , +)$ into the endomorphism ring
of the additive group scheme underlying ${\mathbb L}$ satisfying
the following property: For some trivialization of ${\mathbb L}$
by open affine subschemes ${\rm Spec} \, B$ of $S$ and for each
nonzero $f \in A$ we have $\Phi (f) \mid {\rm Spec} \, B =
{\sum}_{0 \leq i \leq N(f)} \, a_i \cdot \tau^i$ with $a_i \in
B$, such that
\begin{itemize}
\item[(i)] $A \to B$ takes $f$ to $a_0$
\item[(ii)] $a_{N(f)}$ is a unit
\item[(iii)] $N(f) = r \cdot \deg (f)$, where $q^{\deg (f)} = \sharp \,
(A/f)$.
\end{itemize}
To summarize these conditions, we must have $End_{{\mathbb
F}_q}{\mathbb L}=\oplus H^0(S,{\mathbb L}^{q^i -1})\tau$.

A morphism of Drinfeld modules $\Phi \to \Phi'$ over $L$ is an
element $u \in L \{ \tau \}$ such that for all $f \in A$ we have
$u \circ \Phi (f) = \Phi' (f) \circ u$. The endomorphism ring $L
\{ \tau \}$ is the subring $L \{ \tau_p \}$ which is generated by
$\tau : x \mapsto x^q$. Let $n$ be an ideal of $A$. The group
scheme $_n\Phi$ of $n$-division points ${\cap}_{f \in n} \, {\rm
Ker} (\Phi (f))$ is a finite flat subscheme of $({\mathbb L} ,
+)$ of degree $\sharp \, (A/n)^r$ over $S$ which is etale outside
support of $n$. An isomorphism $\alpha : \,_n\Phi \to (A/n)^r$ is
called an $n$-level structure outside the support of $n$.
Equivalently, an $n$-level structure is defined to be a morphism
$\alpha$ from the constant scheme of $A$-modules $(n^{-1} / A)^r$
to $_n\Phi$ such that, ${\sum}_{i \in (n^{-1} / A)^r} \, \alpha
(i) = \, _n\Phi$ as Cartier divisors on ${\mathbb L}$. This
extends the previous definition of an $n$-level structure to
support of $n$.

Let $m$ and $n$ be ideals of $A$ each divisible by at least two
different primes. Such ideals are called admissible ideals. The
moduli scheme over ${\rm Spec} \, A$ which classifies Drinfeld
modules $(\Phi , \alpha)$ of rank $r$ over $S$ together with an
$n$-level structure $\alpha$ will be denoted by $M^r (n)(S)$ or
simply by $M^r (n)$. The moduli space $M^r (n)$ is a smooth affine
scheme of finite type of dimension $r-1$ over $A$. The algebra
$M_n^r$ denotes the affine algebra of the affine scheme $M^r (n)$.
For $n \mid m$ the natural morphism $M^r (m) \to M^r (n)$ is
finite flat and even etale outside the support of $m$. We intend
to work with $M^r (n)$ only over ${\rm Spec} (A [n^{-1}])$. For a
more detailed discussion see [De-Hu] and [Dr1]. We would like to
compactify the etale covers $M^r (n)$ of $M^r (1)$. We require
that for $m \mid n$ there be natural maps between the
compactifications $M^{*r} (m) \to M^{*r} (n)$. At the moment an
analytic theory of Satake compactifications is available [Ge], but
we need an arithmetic theory of compactifications.

In the case of rank-two Drinfeld modules, the moduli space is a
curve and we can compactify it by addition of a few points. This
compactification is functorial. This is the compactification which
Gekeler-Reversat uniformize its Jacobian [Ge-Re].

\subsection{Geometric Drinfeld modular forms}\label{ssec1.2}

Let $R$ be an $A$-algebra. A modular form $f$ of weight $k$ with
respect to ${\rm GL} (r,A)$ is a rule which assigns to each Drinfeld
module $(L,\Phi)$ of rank $r$ over an $R$-scheme $S$ a section
$f(L,\Phi) \in \Gamma (L^{-k})$ with the following property: For any
map of $R$-schemes $g : S' \to S$ the section $f(L,\Phi)$ is
functorial with respect to $g$, {\it i.e.} for any nowhere-zero
section $\beta$ of $g^* (L)$ the element $f(g^* (L) , g^* (\Phi))
\cdot \beta^{\otimes k} \in \Gamma (S' , O_{S'})$ depends only on
the isomorphism class $(E,\Phi,\beta)$. The same definition works
for an arbitrary level structure.

Let $H^r$ denote the graded ring of Drinfeld modular forms of rank
$r$ over $A$, and $H_n^r$ denote the graded ring of Drinfeld modular
forms of rank $r$ and level $n$ over $A$. The map ${\rm Spec}
(H_n^r) \to M^r (n)$ is a principal ${\mathbb G}_m$-bundle which
represents isomorphism classes of Drinfeld modules with level
structures, together with a nowhere-zero section of $L$. Let
$\omega$ denote the top wedge of the sheaf of relative differential
forms on ${\rm Spec} (H_n^r)$ restricted to $M^r (n)$. In rank-two
case, $\omega^{\otimes 2}$ coincides with $\Omega^1 (2 \cdot {\rm
cusps}$).

Here we develop an algebraic theory of ${\mathbf q}$-expansions
for Drinfeld modular forms which uses Tate uniformization (see
[Go]). Let $R$ be a d.v.r. over $A$ with maximal ideal $(\pi)$ and
fraction field $K'$. Let $({\mathbb M} , \Phi)$ be a Drinfeld
module of rank $r$ over $K$. We say that $({\mathbb M} , \Phi)$
has stable reduction modulo $(\pi)$ if for some $c \in K'$ we have
\begin{itemize}
\item[(i)] The module $({\mathbb M} , c \, \Phi \, c^{-1})$ has
coefficients in $R$.
\item[(ii)] The reduction modulo $(\pi)$ is a Drinfeld module of rank
$\leq r$.
\end{itemize}

By definition, equality holds in the good reduction case. A Drinfeld
module with an admissible $n$-level structure, has stable reduction.
Over $A [n^{-1}]$ any $n$-level structure gives us a Drinfeld module
with stable reduction.

\medskip

\noindent {\bf Definition.} {\it Let $K^s$ denote the maximal abelian
extension of $K$ split totally at $\infty$. A $\Phi$-lattice $N$ over
$K$ is a finitely generated projective $A$-submodule of $K^s$ such
that:
\begin{itemize}
\item[(i)] The group $N$ is ${\rm Gal} (K^s / K)$-stable.
\item[(ii)] In any ball, there are only finitely many elements of
$N$. $\Box$
\end{itemize}
}

\medskip

\noindent {\bf Theorem 1.2.1 (Drinfeld).} {\it The isomorphism
classes of Drinfeld modules of rank $r+r_1$ with stable reduction,
are in one-to-one correspondence with the isomorphism classes of
$({\mathbb M} , \Phi , N)$, where $({\mathbb M} , \Phi)$ is a
Drinfeld module of rank $r$ over $K$, with good reduction, and $N$ is
$\Phi$-lattice of rank $r_1$.}

\medskip

One can show that the proper scheme $M^1 (n) \to {\rm Spec} (A)$
is Spectrum of the ring of integers of $K^s$. Let $m$ and $n$ be
ideals of $A$, and $p$ be a prime dividing $m$. Consider the
universal module $\Phi$ over $M^1 (m) / A [p^{-1}]$ with trivial
bundle. For any ring $R$ we define $R(({\mathbf q}))$ to be the
ring of finite tailed Laurent series over $R$. Theorem~1.1
associates a Drinfeld module $T(m,n,p)$ over $(M_m^1 \otimes K)
(({\mathbf q}))$ to the triple $({\mathbb K} , \Phi , \Phi (mn)
(1/ {\mathbf q}))$. The module $T(m,n,p)$ can be extended over
$M_m^1 [p^{-1}] \otimes R (({\mathbf q}))$ with a no-where zero
section ``1'' and a natural $m$-level structure $\psi$. For more
details see [Go]. We call
$$
f(T(m,n,p) , \hbox{``1''} , \psi) \in
H^r \otimes M_m^1 [p^{-1}] (({\mathbf q}))
$$
the ${\mathbf q}$-expansion of $f$ at the cusp $(m,n,p)$. $f$ is
holomorphic at the cusp $(m,n,p)$ if the expansion contains no
negative terms.

\subsection{Analytic Drinfeld modular forms}\label{ssec1.3}

We shall first give an analytic description of Drinfeld modules. We
introduce a construction similar to a ``Weierstrass-preparation''
which starts with introdu\-cing lattices in $C$. An $A$-lattice of rank
$r$ in $C$ is a finitely generated and thus projective $A$-submodule
$\Lambda$ of $C$ of projective rank $r$, whose intersection with each
bounded subset of $C$ is finite. To an $r$-lattice we associate an
exponential function $e_{\Lambda} (z) = z \prod_{\lambda \in \Lambda
- \{ 0 \}} (1-z/\lambda)$. This product converges and defines a
surjective and entire ${\mathbf F}_q$-linear function $e_{\Lambda} :
C \to C$. For $a \in A$ there exists $\Phi_a^{\Lambda} \in C \{ \tau
\}$ such that $e_{\Lambda} (az) = \Phi_a^{\Lambda} (e_{\Lambda}
(z))$. The endomorphism $\Phi_a^{\Lambda}$ induces a $A$-module
structure on the additive group scheme
$$
{\mathbb G}_a / C = C \leftarrow C / {\rm ker} (r_{\Lambda}) = C
/ \Lambda.
$$
In fact $A \mapsto \Phi_a^{\Lambda}$ is a Drinfeld
module. This association is one-to-one. More precisely, the
associated Drinfeld modules are isomorphic if the lattices are
similar.

Let $Y$ be a projective $A$-module of rank $r$. Let $\Gamma_Y$ be the
automorphism group of $Y$. Fixing a $K$-basis for $Y \otimes K$ we
can assume $\Gamma_Y$ is a subgroup of ${\rm GL} (r,K)$ commensurable
with ${\rm GL} (r,A)$. Let $\tilde\Omega^r$ denote $C^r$ with all
hyperplanes defined over $K_{\infty}$ removed. One can see that
$\Gamma_Y \backslash \tilde\Omega^r$ is in one-to-one correspondence
with $A$-lattices isomorphic to $Y$. The quotient $\Omega^r :=
\tilde\Omega^r / C^*$ is a subset of projective $(r-1)$-space. One
can identify $\cup_Y \Gamma_Y \backslash \Omega^r$ with $M^r (1)(C)$.
The union is over a finitely many isomorphism classes of projective
$A$-modules of rank $r$.

In the special case where $A = {\mathbf F}_q [t]$ every finitely
generated projective $A$-module is free. So we can assume
$\Gamma_Y = {\rm GL} (r, A)$ and thus can identify $M^r (1)(C)$
with ${\rm GL} (r, A) \backslash \Omega^r$. In rank two case, let
$\Lambda$ be the 2-lattice $A\omega + A$. We write
$\Phi_T^{\Lambda} = T + g (\omega) \tau + \Delta (\omega)
\tau^2$. The $g (\omega)^{q+1} / \Delta (\omega)$ defines a map
from $\Omega$ to $C$ which gives us a $j$-invariant $j : \Gamma
\backslash \Omega = \Gamma \backslash (C-K_{\infty}) \to M^2
(1)(C) = C$.

Let $k$ be a non-negative integer and let $m$ be an integer class
modulo $q-1$. We define a Drinfeld modular form of rank two, weight
$k$ and type $m$ with respect to the arithmetic group $\Gamma$ to be
a map $f : \Omega \to C$ which is holomorphic in the rigid analytic
sense, such that
$$
f(z\mid \gamma) = (\det \gamma)^{-m} (cz + d)^k f(z) \quad \hbox{for}
\ \gamma = \begin{pmatrix} a &b \\ c &d \end{pmatrix} \in \Gamma \, .
$$
We also require that $f$ be holomorphic at the cusps of $\Gamma$. The
space of these modular forms is denoted by $M_{k,m} (\Gamma)$ and the
subspace of $i$-cuspidal forms is denoted by $M_{k,m}^i (\Gamma)$. We
shall give some cuspidal conditions to determine the exact subspace
which comes from the geometric definition.

For general rank $r$, we define a Drinfeld modular form of weight $k$
and type $Y$ and level $n$ to be a rigid analytic holomorphic
function satisfying the following transformation rule:
$$
f((z,1) \cdot g^t) = f(za + b^t / c^t \cdot z + d) (z \cdot c^t +
d)^{-k} \quad \hbox{for} \ g \in \Gamma_Y (n)
$$
where
$g^t = \begin{pmatrix} a &b \\ c &d \end{pmatrix}$ and $a$ is the
upper left $(r-1) \times (r-1)$-minor. Here we think of $z \in
\Omega^r$ as an element in $A^{r-1}$. The dot product on the left
hand side is the group action.

In case of rank 2 one can give cusp conditions on analytic
Drinfeld modular forms which guarantee algebraicity. Let $f$
denote an analytically defined Drinfeld modular form of weight
$k$ and level $n$. The form $f$ is algebraic if and only if at
each cusp the ${\mathbf q}$-expansion is finite tailed [Go]. If
there are no negative terms at the cusps, then $f$ comes from a
section of $\omega^k$. The space of Drinfeld modular forms
holomorphic at cusps, is finite dimensional.

\subsection{Harmonic cochains}\label{ssec1.4}

First of all, we shall define the Bruhat-Tits tree ${\mathcal I}$ of
${\rm PSL} (2 , K_{\infty})$. The vertices of ${\mathcal I}$ are
defined to be rank two $O_{\infty}$-lattices in $K_{\infty}^2$. For
vertices $L$ and $L'$, we consider the similarity classes $[L]$ and
$[L']$ adjacent, if for some re\-presentatives $L_1 \in [L]$ and $L'_1
\in [L']$ we have $L_1 \subset L'_1$ with $L'_1 / L_1$ being a length
one as $O_{\infty}$-module. The graph ${\mathcal I}$ is a
$(q_{\infty} + 1)$-regular tree. Let $G$ denote the group scheme
${\rm GL} (2)$ with center $Z$. The group $G (K_{\infty})$ acts on
${\mathcal I}$. The stabilizer of the standard vertex
$[O_{\infty}^2]$ is $G(O_{\infty}) \cdot Z (K_{\infty})$. This gives
the standard identifications of the vertices $V({\mathcal I})$ with
$G (K_{\infty}) / G(O_{\infty}) \cdot Z (K_{\infty})$ and the edges
$E({\mathcal I})$ with $G (K_{\infty}) / \Gamma_0 (O_{\infty}) \cdot
Z (K_{\infty})$. Here
$$
\Gamma_0 (O_{\infty}) = \left\{
\begin{pmatrix} a &b \\ c &d \end{pmatrix} \in G (O_{\infty}) \mid c
\equiv 0 \, ({\rm mod} \, \infty) \right\}.
$$

There exists a bijection between the ends of ${\mathcal I}$ and
${\mathbb P}^1 (K_{\infty})$. Also by a theo\-rem of Goldman and
Iwahori [Go-Iw], there exists a $G (K_{\infty})$-equivariant
bijection between ${\mathcal I} ({\mathbb Q})$ and the set of
similarity classes of non-archimedean norms on $K_{\infty}^2$.
Based on this bijection, one can find a $G
(K_{\infty})$-equivariant building map $\lambda : \Omega \to
{\mathcal I} ({\mathbb R})$ with $\lambda (\Omega) = {\mathcal I}
({\mathbb Q})$ (for definition see [Ge-Re] p.~37). By constructing
an explicit open  covering of $\Omega$ (see [Ge-Re] p.~33) one can
define $\tilde\Omega$ and an analytic reduction map $R : \Omega
\to \tilde\Omega$ is a scheme locally of finite type over the
residue field $k_{\infty}$. Each irreducible component of
$\tilde\Omega$ is isomorphic with ${\mathbb P}^1 (k_{\infty})$
and meets exactly $q_{\infty} + 1$ other components. Via the
building map one identifies ${\mathcal I}$ with the intersection
graph of $\tilde\Omega$.

Now we define harmonic cochains. Let $X$ be an abelian group. A
harmonic cochain with values in $X$ is a map $\phi : V({\mathcal I})
\to X$ which satisfies
\begin{itemize}
\item[(i)] $\phi (e) + \phi (\bar e) = 0$ for all $e \in E({\mathcal
I})$ where $\bar e$ denotes the reversely oriented $e$.
\item[(ii)] $\sum_{e \in E ({\mathcal I}) , {\rm tail} (e) = v} \phi
(e) = 0$ for all $v \in V({\mathcal I})$.
\end{itemize}

The group of $X$-valued harmonic cochains is denoted by $\underline H
({\mathcal I} , X)$. Van der Put gives a map $r : O_{\Omega}
(\Omega)^* \to \underline H ({\mathcal I} , {\mathbb Z})$ and obtains
an exact sequence of $G(K_{\infty})$-modules
$$
0 \to C^* \to O_{\Omega} (\Omega)^* \to \underline H ({\mathcal I} ,
{\mathbb Z}) \to 0 \, .
$$
He defines $r$ using the standard covering of $\Omega$ (see
[Ge-Re] p.~40). There is also a residue map from holomorphic
differential forms on $\Omega$ to the space of harmonic cochains
$\Omega^1 (\Omega) \to \underline H ({\mathcal I} , C)$. The
residue map is also defined using the open covering of $\Omega$.
The fact that residue map is well-defined is a consequence of
residue theorem in rigid analytic geometry. The map defined by
Van der Put coincides with the residue map if we reduce these
maps modulo $p$. In other words the following diagram is
commutative (see [Ge-Re] p.~40):
$$
\begin{CD}
O_{\Omega} (\Omega)^* @>{}>> \underline H ({\mathcal I} , {\mathbb
Z}) \\
@V{d \log}VV @VV{\rm res}V \\
\Omega^1 (\Omega) @>{}>> \underline H ({\mathcal I} , C) \, .
\end{CD}
$$

Here $d\log : f \mapsto df / f$ denotes the logarithmic
differentiation. For an arithmetic subgroup $\Gamma$ of $G(K)$, we
define $\underline H ({\mathcal I} , X)^{\Gamma}$ to be the subgroup
of invariants under $\Gamma$. The compact support cohomology
$\underline H_! ({\mathcal I} , X)^{\Gamma}$ is the subgroup of
harmonic cochains which have compact support modulo $\Gamma$. For any
commutative ring $B$ the image of the injective map $\underline H_!
({\mathcal I} , {\mathbb Z})^{\Gamma} \otimes B \to \underline H_!
({\mathcal I} , B)^{\Gamma}$ is denoted by $\underline H_{!!}
({\mathcal I} , B)^{\Gamma}$. In fact $\underline H ({\mathcal I} ,
{\mathbb Z})^{\Gamma}$ is a free abelian group of rank $g$, where $g$
is $\dim_{\mathbb Q} (\Gamma^{ab} \otimes {\mathbb Q})$. Let
$\Gamma^*$ denote $\Gamma$ divided by its torsion subgroup, and
$\tilde\Gamma$ denote the quotient of $\Gamma$ by the finite subgroup
$\Gamma \cap Z(K)$. There exists a map
$$
H_1 (\Gamma \backslash {\mathcal I}) \to \underline H ({\mathcal I} ,
{\mathbb Z})^{\Gamma}
$$
$$
\phi \mapsto \phi^*
$$
defined by $\phi (e) = n(e) \, \phi (\tilde e)$ where $n(e) :=
\sharp \, (\Gamma \cap Z(K))^{-1} \, \sharp \, (\gamma_e)$. We
get an injection with finite cokernel (for definition of these
maps see [Ge-Re] p.~49)
$$
\bar\Gamma = \Gamma^{ab} / {\rm tor} (\Gamma^{ab}) \to
(\Gamma^*)^{ab} \to H_1 (\Gamma \backslash {\mathcal I}) \to
\underline H ({\mathcal I} , {\mathbb Z})^{\Gamma} \, .
$$
Let $c(- , -)$ denote the unique geodesic connecting two vertices.
One can explicitly define the above map. We fix a vertex $v$ of
${\mathcal I}$. For $e \in E({\mathcal I})$ and $\alpha , \gamma \in
\Gamma$ we put $i (e,\alpha , \gamma , v) = 1, -1 , 0$, if $\gamma
(e)$ belongs to $c(v,\alpha (v)), c(\alpha (v),v)$ and neither one,
respectively. Now we get the function
$$
\phi_{\alpha} = \phi_{\alpha , v} := z (\Gamma)^{-1} \sum_{\gamma \in
\Gamma} i(e,\alpha , \gamma , v) \in \underline H_! ({\mathcal I} ,
{\mathbb Z})^{\Gamma}
$$
which is independent of the choice of $v$. One can show that
$\phi_{\alpha , \beta} = \phi_{\alpha} + \phi_{\beta}$. We get an
injection $j : \bar\Gamma \to \underline H_! ({\mathcal I} ,
{\mathbb Z})^{\Gamma}$ with finite cokernel. It is shown that $j$
is an isomorphism in case $K = {\mathbf F}_q (t)$. Also Gekeler
Reversat prove that $j$ is an isomorphism whenever $\tilde\Gamma$
has only $p$-torsion as its torsion subgroup (see [Ge-Re] p.~74).

The point of all this machinery is the following fact (see [Ge-Re]
p.~52, 72):

\medskip

\noindent {\bf Theorem 1.4.1 (Gekeler-Reversat)} {\it Let $\Gamma$ be
an arithmetic subgroup of \break $G(K)$. The map ${\rm res} : f
\mapsto {\rm res} (f)$ induces an isomorphism
$$
M_{2,1}^1 (\Gamma)(C)
\to \underline H_! ({\mathcal I} , C)^{\Gamma}
$$
and an isomorphism
$M_{2,1}^2 (\Gamma) (C) \to \underline H_{!!} ({\mathcal I} ,
C)^{\Gamma}$.}

\medskip

One can identify $M_{2,1}^2 (\Gamma) (C)$ with the vector space
$H^0 (\bar M_{\Gamma} , \Omega^1 ({\rm cusps}))$ and $M_{2,1}^2
(\Gamma)$ $(C)$ with the vector space $H^0 (\bar M_{\Gamma} ,
\Omega^1)$. Formulation of $\underline H ({\mathcal I} ,
B)^{\Gamma}$ and $\underline H_! ({\mathcal I} , B)^{\Gamma}$
commutes with flat ring extensions $B'/B$. Hence we have an
${\mathbf F}_p$-structure $M_{2,1}^2 (\Gamma)$ $({\mathbf F}_p)$
on $M_{2,1}^2 (\Gamma) (C)$.

\subsection{Adelic automorphic forms}\label{ssec1.5}

The basic references for automorphic modular forms are [Ja-La],
[Gel]. Harmonic cochains are related to automorphic forms. To
explain the automorphic theory, we shall first introduce the
Adelic formulation of Drinfeld modular forms. Let ${\mathbb U} =
{\mathbb U}_K = {\mathbb U}_f \times K_{\infty}$ denote the ring
of Adeles over $K$. The Adele group is a locally compact ring
containing $K$ as a discrete cocompact subring.

The ring of integers ${\mathbb O} = {\mathbb O}_K = {\mathbb O}_f
\times O_{\infty}$, is the maximal compact subring of ${\mathbb U}$.
We have a decomposition of idle group ${\mathbb I} = {\mathbb I}_K =
{\mathbb I}_f \times K_{\infty}^*$ where ${\mathbb I}_f$ denotes the
finite idles. The ring ${\mathbb O}_f$ is the completion with respect
to ideal topology of $A = K \cap {\mathbb O}_f$. The class group of
$A$ can be identifies with the quotient $K^* \backslash {\mathbb I}_f
/ {\mathbb O}_f^*$. Hence this quotient is finite. We have a
bijection between $K^* \backslash {\mathbb I}_f / {\mathbb O}_f^*$
and $G(K) \backslash G({\mathbb U}_f) / G({\mathbb O}_f)$ where $G =
{\rm GL} (2)$. The latter is identified with isomorphism classes of
rank-two $A$-lattices in $K^2$ in the following manner. We identify
$g \in G ({\mathbb U}_f)$ with the $A$-lattice whose span in
${\mathbb O}_f^2$ is the same as $g^{-1} \cdot {\mathbb O}_f^2$. Here
the action of $g$ is action of a matrix. More generally as a result
of strong approximation theorem for the group ${\rm SL} (2)$, for any
open subgroup ${\mathbb K}_f$ in $G({\mathbb U}_f)$ determinant
induces a bijection
$$
G(K) \backslash G({\mathbb U}_f) / {\mathbb K}_f \to K^* \backslash
{\mathbb I}_f / \det ({\mathbb K}_f) \, .
$$
Choosing a measure $\mu$ on the locally compact group $G({\mathbb U})
/ Z({\mathbb U})$ we always get a finite volume for $G(K) \backslash
G({\mathbb U}) / Z({\mathbb U})$. So is the case for $G(K) \backslash
G({\mathbb U}) / Z(K_{\infty})$.

Now consider the Hilbert space of complex valued square integrable
functions on $G(K) \backslash G({\mathbb U}) / Z(K_{\infty})$ and
denote it by $L^2 (G(K) \backslash G({\mathbb U}) / Z(K_{\infty}))$.
The subspace $L_!^2 (G(K) \backslash G({\mathbb U}) / Z(K_{\infty}))$
consists of $G({\mathbb U})$-stable functions $\phi$ which satisfy
the cusp condition
$$
\int_{K \backslash {\mathbb U}} \Phi \left( \begin{pmatrix} 1 &x \\ 0
&1 \end{pmatrix} \cdot g \right) dx = 0 \, .
$$
Here $dx$ is a Haar measure on $K \backslash {\mathbb U}$. We have a
decomposition $L_!^2 = \hat\oplus \, L_!^2 (\psi)$, according to the
character $\psi$ of the compact group
$$
Z(K) \backslash Z({\mathbb U}) / Z(K_{\infty}) = K^* \backslash
{\mathbb I} / K_{\infty}^*.
$$
For each $\psi$ we have $L_!^2
(\psi) = \hat\oplus \, V_Q$, where $V_Q$'s are irreducible unitary
$G({\mathbb U})$-submodules occuring with multiplicity one. The
underlying unitary representation will be called cuspidal
automorphic representations. In each $V_Q$ there exists a
distinguished nonzero function $\phi_Q$, satisfying $\phi (gk) =
\phi (g) \, \psi (k)$ for all $g \in G({\mathbb U})$ and $k \in
{\mathbb K}_Q (n_Q)$. The open subgroup ${\mathbb K}_Q (n_Q)
\subset G({\mathbb O})$ is defined by
$$
{\mathbb K}_Q (n_Q) = \left\{ \begin{pmatrix} a &b \\ c &d
\end{pmatrix} \in G({\mathbb O}) \mid c \equiv 0 \, ({\rm mod} \,
n_Q) \right\}
$$
where $n_Q$ is an integer, the conductor of the representation. Here
$\psi$ is extended to ${\mathbb K}_Q (n_Q)$ by defining $\psi (k) =
\prod_v \psi_v (a_v)$, where $k = \begin{pmatrix} a &b \\ c &d
\end{pmatrix}$ and $a=(a_v)$.

The support of automorphic cusp forms lives on the double coset
space $G(K) \backslash$ $G({\mathbb U}) / {\mathbb K}_Q (n_Q)
\cdot Z(K_{\infty})$. So the space $W ({\mathbb K}_Q (n_Q))$ of
automorphic cusp forms for ${\mathbb K}_Q (n_Q)$ has finite
dimension and is independent of ${\mathbb C}$, {\it i.e.} for any
field $F$ of characteristic zero $W({\mathbb K}_Q (n_Q) , F)
\otimes {\mathbb C} = W ({\mathbb K}_Q (n_Q) , {\mathbb C})$.

Only a subspace of $W ({\mathbb K}_f^0 (n))$ is related to harmonic
cochains. Consider the subspace of representations transforming like
the irreducible $G(K_{\infty})$-module
$$
V_{sp,F} = \{ f : {\mathbb P}^1 (K_{\infty}) \to F \mid f \
\hbox{locally constant}\} / F \, .
$$
The underlying representation is independent of $F$ and is denoted by
$Q_{sp}$. Over ${\mathbb C}$ there is a unique (up to scaling)
$G(K_{\infty})$-invariant scalar product on $V_{sp}$ such that
$Q_{sp}$ extends to $\hat Q_{sp}$ on $\hat V_{sp}$. The conductor is
the prime divisor $\infty$. An $F$-valued automorphic form is said to
transform like $Q_{sp}$ if the space of its $G(K_{\infty})$-right
translates generate a module isomorphic to a finite number of copies
of $Q_{sp}$.

\medskip

\noindent {\bf Theorem 1.5.1 (Drinfeld) } [Dr1] {\it For an open
subgroup ${\mathbb K}_f$ of $G({\mathbb O}_f)$ and for ${\mathbb
K} = {\mathbb K}_f \times \Gamma_0 (\infty)$ there exists a
canonical bijection
$$
\oplus_{x\in S} \, \underline H_! ({\mathcal
I} , F)^{\Gamma_x} \to W_{sp} ({\mathbb K} , F)
$$
where $S$ is a choice of representatives for the finite set $G(K)
\backslash G({\mathbb U}_f) / {\mathbb K}_f$.}

\subsection{Hecke operators}\label{ssec1.6}

In this section we give an analytic description of Hecke operators
and an algebraic description of Hecke correspondences. Again this
picture is motivated by the traditional adelic definition of Hecke
correspondences using double cosets. Let ${\mathbb K}_f$ be an open
subset in $G({\mathbb O}_f)$ with conductor $n$. The conductor is the
least positive divisor coprime to $\infty$ such that ${\mathbb K}_f
(n) \subset {\mathbb K}_f$. Let $v$ be a finite place coprime to $n$.
Then ${\mathbb K}_v := G(O_v)$ embeds in ${\mathbb K}_f$. If $\pi_v$
is a local uniformizer, and $\tau_v = {\rm diag} \, (\pi_v , 1) \in
G(K_v)$, the group ${\mathbb H}_v = {\mathbb K}_v \cap \tau_v \,
{\mathbb K}_v \, \tau_v^{-1}$ has index $q_v + 1$ in ${\mathbb K}_v$.
We define the Hecke operator acting on a function $\phi$ on $G(K)
\backslash G({\mathbb U}) / {\mathbb K} \cdot Z (K_{\infty})$ by the
integral
$$
(T_v \, \phi)(g) := \int_{{\mathbb K}_v} \phi (g \, k_v \, \tau_v) \,
dk_v = \sum_{k_v \in {\mathbb K}_v \backslash {\mathbb H}_v} \phi (g
\, k_v \, \tau_v) \, .
$$
Here $k_v$ runs through a set of representatives of ${\mathbb K}_v
\backslash {\mathbb H}_v$ and $dk_v$ is a Haar measure chosen in a
way that the volume of the quotient ${\mathbb K}_v \backslash
{\mathbb H}_v$ be equal to one.

The Hecke operators preserve $W_{sp}$. They commute. Together with
their adjoints generate a commutative algebra of normal operators
on $W_{sp}$. So there exists a basis of simultaneous eigen-forms
for all $T_v$. The eigenvalues are algebraic integers. Drinfeld
used Weil conjectures to prove analogue of Ramanujan's
conjecture, bounding the norms of eigenvalues $\vert \lambda
(\phi , v) \vert \leq 2 \, q_v^{1/2}$. The Hecke operators
respect the decomposition $W_{sp} = \oplus_{\psi} \, W_{sp}
(\psi)$. We know that
$$
T_v^* \mid W_{sp}(\psi) = \psi^{-1} (\pi_v) \cdot T_v \mid W_{sp}.
$$

Let $K_{\infty}^{\rm un}$ be the maximal unramified extension of
$K_{\infty}$. The Galois group ${\rm Gal} (K_{\infty}^{\rm un} /$
$K_{\infty})$ is isomorphic to the profinite completion
$\hat{\mathbb Z}$, where the canonical generator corresponds to
the Frobenius element ${\rm Fr}_{\infty}$ of the extension
$K_{\infty}^{\rm un}$ over $K_{\infty}$. Let $l$ be a prime
different from $p = {\rm char} ({\mathbf F}_q)$ and let $E_l /
K_{\infty}^{\rm un}$ be the extension obtained by adding all the
$l^r$-th roots of the uniformizer $\pi_{\infty}$ to
$K_{\infty}^{\rm un}$. Then ${\rm Gal} (E_l / K_{\infty}^{\rm
un}) = {\mathbb Z}_l (1)$. Moreover $E_l$ is Galois over
$K_{\infty}$ and we have
$$
{\rm Gal} (E_l / K_{\infty}) = {\rm Gal} (K_{\infty}^{\rm un} /
K_{\infty}) \propto {\rm Gal} (E_l / K_{\infty}^{\rm un}) =
\hat{\mathbb Z} \propto {\mathbb Z}_l (1) \, .
$$
For any $u \in {\mathbb Z}_l$ the action of ${\rm Fr}_{\infty} = 1
\in \hat{\mathbb Z}$ on ${\mathbb Z}_l$ is given by ${\rm
Fr}_{\infty} \, u \, {\rm Fr}_{\infty} = u^{q_{\infty}}$. We get a
two dimensional Galois representation
$$
sp = sp_l = {\rm Gal} (K_{\infty}^{\rm sep} / K_{\infty}) \to {\rm
Gal} (E_l / K_{\infty}) = \hat{\mathbb Z} \propto {\mathbb Z}_l \to
{\rm GL} (2 , {\mathbb Q}_l) \, .
$$
Here we have chosen an isomorphism between ${\mathbb Z}_l (1)$ and
${\mathbb Z}_l$. The last arrow maps $(1,0)$ to $\begin{pmatrix} 1 &0
\\ 0 &q_{\infty}^{-1} \end{pmatrix}$ and $(0,1)$ to $\begin{pmatrix}
1 &1 \\ 0 &1 \end{pmatrix}$.

Now we can state the Drinfeld reciprocity law [Ge2]:

\medskip

\noindent {\bf Theorem 1.6.1 (Drinfeld).} {\it Let $\bar M_{{\mathbb
K}_f}$ denote the compactified Drinfeld moduli space over $C$
associated to ${\mathbb K} = {\mathbb K}_f \times {\mathbb
I}_{\infty}$. There exists a Hecke invariant isomorphism between
$W_{\rm sp} ({\mathbb K} , {\mathbb Q}_l) \otimes {\rm sp}_l$ and $H^1
(\bar M_{{\mathbb K}_f} , {\mathbb Q}_l)$ which is compatible with the
local Galois action of ${\rm Gal} (K_{\infty}^{\rm sep} /
K_{\infty})$.}

\medskip

\subsection{Theta functions}\label{ssec1.7}

A holomorphic theta function $u : \Omega \to C$ with respect to
$\Gamma$ is a holomorphic nonzero function on $\Omega$ and on cusps,
such that $u(\alpha z) = c_u (\alpha) \, u(z)$ for some $c_u (\alpha)
\in C^*$ and all $\alpha \in C$. For a meromorphic theta function, we
allow poles and zeroes on $\Omega$ but not at cusps. The logarithmic
derivative of a holomorphic theta function is in $M_{2,1}^2 (\Gamma)$.
Let $\Theta_m (\Gamma)$ denotes the multiplicative group of meromorphic
theta functions and $\Theta_h (\Gamma)$ denotes the subspace of
holomorphic theta functions. The family of maps $c_u : \Gamma \to C^*$
induces a map
$$
c : \Theta_m (\Gamma) \to {\rm Hom} (\Gamma , C^*) = {\rm Hom}
(\Gamma^{ab} \to C^*).
$$
We have $\ker (c) \cap \Theta_h (\Gamma) = C^*$. For fixed
$\omega , \eta \in \Omega$ define the following holomorphic theta
function
$$
\theta_{\Gamma} (\omega , \eta , z) = \prod_{\gamma \in \tilde\Gamma}
\frac{z - \gamma \omega}{z - \gamma \eta} \, .
$$
Here $\tilde\Gamma = \Gamma/\Gamma \cap Z(K)$. For $\alpha \in \Gamma$
the constant $c(\omega , \eta , \alpha) \in C^*$ defined by
$$
\theta_{\Gamma} (\omega , \eta , \alpha z) = c (\omega , \eta , \alpha)
\cdot \theta_{\Gamma} (\omega , \eta , z)
$$
introduces a group homomorphism $\Gamma \to C^*$ with factors
through the quotient $\bar\Gamma = \Gamma^{ab} /$ torsion. Now
consider the holomorphic function
$$
u_{\alpha} (z) = \theta (\omega , \alpha \omega , z) \, .
$$
One can see that this is a holomorphic function independent of $\omega$
and depending only on class of $\alpha$ in $\bar\Gamma$. We also have
$u_{\alpha\beta} (z) = u_{\alpha} (z) \cdot u_{\beta} (z)$. Now the
equality
$$
c (\omega , \eta , \alpha) = u_{\alpha} (\eta) / u_{\alpha} (\omega)
$$
helps us to define a symmetric bilinear map $\bar\Gamma \times
\bar\Gamma \to C^*$ via $(\alpha , \beta) \mapsto c_{\alpha}
(\beta)$ where $c_{\alpha} (\beta) := c(\omega , \alpha \omega ,
\beta) \in K_{\infty}$. See [Ge-Re], 65 for more details.

On the other hand the residue map we defined before $r :
O_{\Omega} (\Omega)^* \to \underline H ({\mathcal I} , {\mathbb
Z})$ satisfies $r(u_{\alpha}) = \phi_{\alpha} \in \underline H_!
({\mathcal I} , {\mathbb Z})^{\Gamma}$. The Peterson pairing on
$\underline H_! ({\mathcal I} , {\mathbb Z})^{\Gamma}$ is
compatible with pairing on $\bar\Gamma$ that we defined above.
More precisely for $\alpha , \beta \in \Gamma$ we have $2(\alpha
, \beta) = (\phi_{\alpha} , \phi_{\beta})_{\mu}$. We get
injectivity of $\bar c : \bar\Gamma \to {\rm Hom} (\bar\Gamma ,
C^*)$ induced by $\alpha \to c_{\alpha}$.

\section{Congruences between Drinfeld modular forms}\label{sec2}

In this paper we prove an Ihara lemma which is a kind of cokernel
torsion-freeness result for Drinfeld modular forms. This can be
used to obtain congruences between Drinfeld modular forms. The
idea of using Ihara lemma to get congruences is due to Ribet [Ri].

\subsection{Uniformization of the Jacobian}\label{ssec2.1}

Gekeler-Reversat in their important paper [Ge-Re] introduce a
uniformization of the Jacobian of Drinfeld modular curve ${\rm
Jac} (\bar M_{\Gamma})$. Here we summarize the procedure of
constructing this uniformization. More details can be found in
section 7 of [Ge-Re]. A subgroup $\Lambda$ of $(C^*)^g$ is called
a lattice, if $\Lambda \simeq {\mathbb Z}^g$ and the image of
$\Lambda$ under $\log : (C^*)^g \to {\mathbb R}^g$ is a lattice
in ${\mathbb R}^g$. The logarithm depends on the choice of a
basis for the character group of $(C^*)^g$, but the lattice
property is independent of this choice. Thus, we can define
lattices in tori over arbitrary complete fields in particular
over $K_{\infty}$. For any torus $T$ over $C$ of dimension $g$
and any lattice $\Lambda$ in $T$, $T/\Lambda$ is an analytic group
variety which is compact in the rigid analytic sense. $T/\Lambda$
is projective algebraic if and only if there exists a
homomorphism from $\Lambda$ to the character group ${\rm Hom} (T
, {\mathbb G}_m)$ which induces a symmetric positive definite
pairing $\Lambda \times \Lambda \to C^*$. Positive definiteness
means that $\log \vert (\alpha , \alpha) \vert > 0$, for $\alpha
\in \Lambda$ which is not equal to 1. The requirement that $C$ is
algebraically closed is not necessary to get algebraicity of
$T/\Lambda$. For any symmetric positive definite pairing $\Lambda
\times \Lambda \to C_0^*$ with $C_0$ any complete subextension of
$C/K_{\infty}$ we get a uniformization of an abelian variety
defined over $C_0$. Given an arithmetic subgroup $\Gamma \subset
{\rm GL} (2,K)$, we have defined a pairing (see \ref{ssec1.7})
$\bar\Gamma \times \bar\Gamma \to C^*$ which induces an analytic
uniformization of some abelian variety $A_{\Gamma} (C)$. We get an
exact sequence
$$
0 \to \bar\Gamma \overset{\bar c}{\to} {\rm Hom} (\bar\Gamma , C^*)
\to A_{\Gamma} (C) \to 0 \, .
$$
Now fix $\omega_0 \in \Omega$ and consider the map $\psi : \Omega
\to {\rm Hom} (\bar\Gamma , C^*) \to A_{\Gamma} (C)$ which is
induced by $\omega \mapsto c (\omega_0 , \omega , \cdot)$. The
map $\psi$ factors through $\Omega \to \Gamma \backslash \Omega =
M_{\Gamma} (C)$ and extends to an analytic embedding $\bar
M_{\Gamma} (C) \to A_{\Gamma} (C)$. By GAGA theorems this is a
morphism of algebraic varieties. Using theta series one can show
that the abelian variety $A_{\Gamma} (C)$ is
$K_{\infty}$-isomorphic to the Jacobian of the compactified
modular curve ${\rm Jac} (\bar M_{\Gamma})$ (see [Ge-Re] p. 77).
So in fact, we have obtained an analytic uniformization of the
Jacobian which has a purely function field origin. The map $\bar
M_{\Gamma} (C) \to A_{\Gamma} (C)$ factors through ${\rm Jac}
(\bar M_{\Gamma})$ and we use the induced isomorphism ${\rm Jac}
(\bar M_{\Gamma}) \overset{\simeq}{\to} A_{\Gamma} (C)$ identify
${\rm Jac} (\bar M_{\Gamma})$ and $A_{\Gamma} (C)$ from now on.
Let $n$ be a divisor in ${\mathbf F}_q [t]$. We define the
following congruence groups:
$$
\Gamma (n) = \left\{ \gamma \in {\rm SL} (2 , {\mathbf F}_q [t])
\mid \gamma \equiv \begin{pmatrix} 1 &0 \\ 0 &1 \end{pmatrix}
\quad ({\rm mod} \, n) \right\}
$$
$$
\Gamma_0 (n) = \left\{ \gamma \in {\rm SL} (2 , {\mathbf F}_q [t])
\mid \gamma \equiv \begin{pmatrix} * &* \\ 0 &* \end{pmatrix}
\quad ({\rm mod} \, n) \right\}
$$
$$
\Gamma^0 (n) = \left\{ \gamma \in {\rm SL} (2 , {\mathbf F}_q [t])
\mid \gamma \equiv \begin{pmatrix} * &0 \\ * &* \end{pmatrix}
\quad ({\rm mod} \, n) \right\} \, .
$$
Let $v$ be a point outside the divisor $n$. Conjugation by
$\tau_v = \begin{pmatrix} \pi_v &0 \\ 0 &1 \end{pmatrix}$ gives an
isomorphism $\overline{\Gamma_0 (v) \cap \Gamma (n)} \to
\overline{\Gamma^0 (v) \cap \Gamma (n)}$. Here $\pi_v$ is a
uniformizer for the prime ideal $v$. Conjugation of $\Gamma_0 (v)$
by $\begin{pmatrix} 0 &1 \\ -1 &0 \end{pmatrix}$ gives another
isomorphism. Composing one with the inverse of the other induces
an involution $w_v$ of the space $\overline{\Gamma_0 (v) \cap
\Gamma (n)}$ which is called the Atkin-Lehner involution. This
induces an involution $w'_v$ of ${\rm Hom} (\overline{\Gamma_0
(v) \cap \Gamma (n)} , C^*)$ and an involution $w''_v$ of ${\rm
Jac} (\bar M_{\Gamma_0 (v) \cap \Gamma (n)})$.

Let $I : \overline{\Gamma_0 (v) \cap \Gamma (n)} \to
\overline{\Gamma (n)}$ denote the map induced by inclusion. The
transfer map has the reverse direction $V : \overline{\Gamma (n)}
\to \overline{\Gamma_0 (v) \cap \Gamma (n)}$. We have $I \circ V
= [ q_v + 1]$. We define the Hecke operator $T_v :
\overline{\Gamma (n)} \to  \overline{\Gamma (n)}$ to be the
composition $I \circ w_v \circ V$. The operator $T_v$ acting on
$ \overline{\Gamma (n)}$ is self adjoint with respect to the
pairing $(\alpha , \beta) \mapsto c_{\alpha} (\beta)$. The action
of $T_v$ commutes with the isomorphisms $J_{\Gamma (n)} :
\overline{\Gamma (n)} \to \underline H_! ({\mathcal I} ,
{\mathbb Z})^{\Gamma (n)}$. The operator $T_v = I \circ w_v \circ V$
induces an action $T_v : {\rm Hom} (\Gamma (n) , C^*) \to {\rm Hom}
(\Gamma (n) , C^*)$.

The Jacobian ${\rm Jac} (\bar M_{\Gamma (n)})$ is also equipped
with action of the Hecke operator $T_v$. Note that $\bar M_{\Gamma
(n)}$ has several geometric components and $T_v$ does not respect
them. Gekeler-Reversat prove that, this Hecke action is compatible
with Hecke actions on $\overline{\Gamma (n)}$ and ${\rm Hom}
(\overline{\Gamma (n)} , C^*)$. So uniformization of the Jacobian
is Hecke equivariant (see [Ge-Re] p. 86). In fact, we have the
following commutative diagram with exact rows, where the right
hand vertical maps are induced by $V$ and $I$:
$$
\begin{CD}
{\scriptstyle 0} @>{}>> {\scriptstyle \overline{\Gamma (n)}}
@>{\bar c}>> {\scriptstyle {\rm Hom} (\overline{\Gamma (n)} , C^*)}
@>{}>> {\scriptstyle {\rm Jac} (\bar M_{\Gamma (n)})(C)}
@>{}>> {\scriptstyle 0} \\
&&@A{I}AA @A{V}AA @A{}AA \\
{\scriptstyle 0} \!\!\!\! @>{}>> \!\!\! {\scriptstyle
\overline{\Gamma_0 (v) \cap \Gamma (n)}} \!\!\! @>{\bar c}>> \!\!\!
{\scriptstyle {\rm Hom} (\overline{\Gamma_0 (v) \cap \Gamma (n)} , C^*)}
\!\!\! @>{}>> \!\!\! {\scriptstyle {\rm Jac} (\bar M_{\Gamma_0 (v)
\cap \Gamma (n)})(C)} \!\!\! @>{}>> \!\!\!\! {\scriptstyle 0} \\
&&@A{V}AA @A{I}AA @A{}AA \\
{\scriptstyle 0} @>{}>> {\scriptstyle \overline{\Gamma (n)}}
@>{\bar c}>> {\scriptstyle {\rm Hom} (\overline{\Gamma (n)} , C^*)}
@>{}>> {\scriptstyle {\rm Jac} (\bar M_{\Gamma (n)})(C)}
@>{}>> {\scriptstyle 0} \, .
\end{CD}
$$
This is because $c_{I(\alpha)} (\beta) = c_{\alpha} (V(\beta))$
and $c_{V(\alpha)} (\beta) = c_{\alpha} (I(\beta))$ (follows from
[Ge-Re] 6.3.2).

\medskip

\noindent {\bf Proposition 2.1.1.} {\it The map ${\rm Jac}
(\bar M_{\Gamma_0 (v) \cap \Gamma (n)})(C) \to {\rm Jac}
(\bar M_{\Gamma (n)}) (C)$ and the map ${\rm Jac}
(\bar M_{\Gamma (n)}) (C) \to {\rm Jac} (\bar M_{\Gamma_0 (v)
\cap \Gamma (n)})(C)$ in the above diagram are the same as the
maps given by Albanese functoriality and Picard functoriality
respectively.}

\medskip

\noindent {\bf Proof.} Fixing $\omega_0 \in \Omega$ and using
Gekeler-Reversat uniformization we get embeddings $\bar
M_{\Gamma_0 (v) \cap \Gamma (n)} \to {\rm Jac} (\bar M_{\Gamma_0
(v) \cap \Gamma (n)})$ and $\bar M_{\Gamma (n)} \to {\rm Jac}
(\bar M_{\Gamma (n)})$. We have to check that the map between
Jacobians ${\rm Jac} (\bar M_{\Gamma_0 (v) \cap \Gamma (n)})(C)
\to  {\rm Jac} (\bar M_{\Gamma (n)}) (C)$ induced by the above
diagram restricts to a projection $\bar M_{\Gamma_0 (v) \cap
\Gamma (n)} \to \bar M_{\Gamma (n)}$ between the embedded curves.
Then automatically the projection will be the Albanese map. The
fact that $C_{I\alpha} (\beta) = c_{\alpha} (V\beta)$ implies
that $\bar M_{\Gamma_0 (v) \cap \Gamma (n)}$ projects to $\bar
M_{\Gamma (n)}$. (For details see the proof of theorem 7.4.1 in
[Ge-Re]). To show that ${\rm Jac} (\bar M_{\Gamma (n)}) (C) \to
{\rm Jac} (\bar M_{\Gamma_0 (v) \cap \Gamma (n)})(C)$ is given by
Picard functoriality we should see if a divisor on $\bar
M_{\Gamma (n)}$ pulls back to the appropriate set of points on
$\bar M_{\Gamma_0 (v) \cap \Gamma (n)}$. We should check that
$\prod c(\omega_0 , \alpha_i \, \omega , \beta) = c (\omega_0 ,
\omega ,V\beta)$, where $\alpha_i$ are representatives of
$\Gamma_0 (v) \cap \Gamma (n) \backslash \Gamma (n)$. This is
proven in proposition 6.3.2 of [Ge-Re]. $\Box$

\bigskip

We have the following commutative diagrams (see [Ge-Re] p. 72):
$$
\begin{CD}
{\scriptstyle \underline H_! ({\mathcal I} , {\mathbb Z})^{\Gamma (n)}}
@<{}<< {\scriptstyle \overline{\Gamma (n)}} &{\scriptstyle
\overline{\Gamma (n)}} @>{}>> {\scriptstyle \underline H_!
({\mathcal I} , {\mathbb Z})^{\Gamma (n)}} \\
@V{}VV @VV{V}V  \!\!\!\!\!\!\!\!\!\!\!\!\!\!\!\!\!\!\!\!\!\!\!\!\!\!\!\!
\!\!\!\!\!\!\!\!\!\!\!\!\!\!\!\!\!\!\!\!\!\!\!\!\!\!\!\!@A{I}AA
\!\!\!\!\!\!\!\!\!\!\!\!\!\!\!\!\!\!\!\!\!\!\!\!\!\!\!\!@AA{\rm Norm}A \\
{\scriptstyle \underline H_! ({\mathcal I} , {\mathbb Z})^{\Gamma_0 (v)
\cap \Gamma (n)}} @<{}<< {\scriptstyle \overline{\Gamma_0 (v) \cap
\Gamma (n)}}  & \ \ {\scriptstyle \overline{\Gamma_0 (v) \cap \Gamma (n)}}
@>{}>> {\scriptstyle \underline H_! ({\mathcal I} ,
{\mathbb Z})^{\Gamma_0 (v) \cap \Gamma (n)}} \, .
\end{CD}
$$

Assume $X - \{ \infty \}$ is the affine line over ${\mathbf
F}_q$. By [Ge-Re] p. 74 there is an isomorphism $\bar\Gamma
\overset{\simeq}{\to} \underline H_! ({\mathcal I} , {\mathbb
Z})^{\Gamma}$. By commutativity of the above diagrams, these
isomorphisms are Hecke-equivariant. The non-degenerate pairing on
$\bar\Gamma$ is compatible with the Petersson pairing on
$\underline H_! ({\mathcal I} , {\mathbb Z})^{\Gamma}$. Also from
the injectivity of $\underline H_! ({\mathcal I} , {\mathbb
Z})^{\Gamma (n)} \to \underline H_! ({\mathcal I} , {\mathbb
Z})^{\Gamma_0 (v) \cap \Gamma (n)}$ we get the injectivity of $V
: \overline{\Gamma (n)} \to \overline{\Gamma_0 (v) \cap \Gamma
(n)}$.

\subsection{${\rm SL} (2 , {\mathbf F}_q [t])$ and congruences}\label{ssec2.2}

In this section we assume that $q > 4$ and $X - \{ \infty \}$ is
the affine line over ${\mathbf F}_q$. Let $\Gamma$ denote an
arbitrary arithmetic subgroup of ${\rm GL}
(2 , {\mathbf F}_q (t))$ and $\bar\Gamma$ denote $\Gamma^{ab}$
divided by its torsion subgroup. We start with the analytic
uniformization of ${\rm Jac} (\bar M_{\Gamma})$ given by the
exact sequence
$$
0 \to \bar\Gamma \overset{\bar c}{\to} {\rm Hom} (\bar\Gamma ,
C^*) \to {\rm Jac} (\bar M_{\Gamma}) (C) \to 0.
$$

Let $\pi : \bar M_{\Gamma_0 (v) \cap \Gamma (n)} \to \bar
M_{\Gamma (n)}$ denote the natural projection induced by the
injection of the corresponding congruence groups. Conjugation by
the matrix $\begin{pmatrix} 0 &1 \\ -1 &0 \end{pmatrix} \cdot
\tau_v$ induces an involution $w_v$ of $\bar M_{\Gamma_0 (v)
\cap \Gamma (n)}$. We have the following commutative diagram
$$
\begin{CD}
{\scriptstyle 0} \!\! @>{}>> \!\! {\scriptstyle
\overline{\Gamma_0 (v) \cap \Gamma (n)}} \!\! @>{\bar c}>> \!\!
{\scriptstyle {\rm Hom} (\overline{\Gamma_0 (v) \cap \Gamma (n)}
, C^*)} \!\! @>{}>> \!\! {\scriptstyle {\rm Jac} (\bar M_{\Gamma_0
(v) \cap \Gamma (n)})(C)} \!\! @>{}>> \!\! {\scriptstyle 0} \\
&&@V{w_v}VV @V{w'_v}VV @V{w''_v}VV \\
{\scriptstyle 0} \!\! @>{}>> \!\! {\scriptstyle \overline{\Gamma_0 (v)
\cap \Gamma (n)}} \!\! @>{\bar c}>> \!\! {\scriptstyle {\rm Hom}
(\overline{\Gamma_0 (v) \cap \Gamma (n)} , C^*)} \!\! @>{}>> \!\!
{\scriptstyle {\rm Jac} (\bar M_{\Gamma_0 (v) \cap \Gamma (n)})(C)}
\!\! @>{}>> \!\! {\scriptstyle 0} \, .
\end{CD}
$$

We define a map $\pi' : \pi \circ w_v : \bar M_{\Gamma_0 (v) \cap
\Gamma (n)} \to  \bar M_{\Gamma (n)}$. We get the following
commutative diagram with exact rows:
$$
\begin{CD}
{\scriptstyle 0} \!\! @>{}>> \!\! {\scriptstyle \overline{\Gamma
(n)}^{\oplus 2}} \!\! @>{\bar c}>> \!\! {\scriptstyle {\rm Hom}
(\overline{\Gamma (n)} , C^*)^{\oplus 2}} \!\! @>{}>> \!\!
{\scriptstyle {\rm Jac} (\bar M_{\Gamma (n)})(C)^{\oplus 2}}
\!\! @>{}>> \!\! {\scriptstyle 0} \\
&&@A{(I,I \circ w_v)}AA @A{(V,V \circ w'_v)}AA @AA{(\pi_* ,
\pi_* \circ w''_v)}A \\
{\scriptstyle 0} \!\! @>{}>> \!\! {\scriptstyle \overline{\Gamma_0
(v) \cap \Gamma (n)}} \!\! @>{\bar c}>> \!\! {\scriptstyle {\rm Hom}
(\overline{\Gamma_0 (v) \cap \Gamma (n)} , C^*)} \!\! @>{}>> \!\!
{\scriptstyle {\rm Jac} (\bar M_{\Gamma_0 (v) \cap \Gamma (n)})(C)}
\!\! @>{}>> \!\! {\scriptstyle 0} \\
&&@A{V \circ \pi_1 + w_v \circ V \circ \pi_2}AA @A{I \circ \pi_1 +
w'_v \circ I \circ \pi_2}AA @AA{\pi^* \circ \pi_1 + w''_v \circ
\pi^* \circ \pi_2}A \\
{\scriptstyle 0} \!\! @>{}>> \!\! {\scriptstyle \overline{\Gamma
(n)}^{\oplus 2}} \!\! @>{\bar c}>> \!\! {\scriptstyle {\rm Hom}
(\overline{\Gamma (n)} , C^*)^{\oplus 2}} \!\! @>{}>> \!\!
{\scriptstyle {\rm Jac} (\bar M_{\Gamma (n)})(C)^{\oplus 2}}
\!\! @>{}>> \!\! {\scriptstyle 0} \, .
\end{CD}
$$

Here $\pi_1$ and $\pi_2$ denote the projections to the first and
second component respectively. This diagram gives us control on the
maps between Jacobians induced by the involution $w''_v$ and
Albanese functoriality or Picard functoriality.

\medskip

\noindent {\bf Theorem 2.2.1.} {\it Let $n$ be an effective divisor
on $X - \{\infty \}$ and let $v$ be a point on $X - \{\infty \}$
which does not intersect $n$. For $l$ not dividing $2q$ the map
$$
\alpha = (I,I \circ w_v) : \overline{\Gamma_0 (v) \cap \Gamma (n)}
\to \overline{ \Gamma (n)} \oplus  \overline{ \Gamma (n)}
$$
is a surjection modulo $l$.}

\medskip

\noindent {\bf Proof.} Serre's machinery in geometric group
theory discusses action of groups on trees [Se]. It can be used to
understand the algebraic nature of $\alpha$. Serre proves that
$$
{\rm SL} (2 , {\mathbf F}_q (t)_v) = {\rm SL} (2 , {\mathbf F}_q
[t]_v) *_{\Gamma^v} {\rm SL} (2 , {\mathbf F}_q [t]_v).
$$
The congruence group $\Gamma^v$ is defined by
$$
\Gamma^v = \left\{
\gamma \in {\rm SL} (2 , {\mathbf F}_q [t]_v) \mid \gamma \equiv
\begin{pmatrix} * &* \\ 0 &* \end{pmatrix} ({\rm mod} \, \pi_v)
\right\}.
$$
Density of ${\mathbf F}_q [t] [v^{-1}]$ in $ {\mathbf
F}_q (t)$ implies that
$$
{\rm SL} (2 , {\mathbf F}_q [t] [v^{-1}]) = {\rm SL} (2 ,
{\mathbf F}_q [t]) *_{\Gamma_0 (v)} {\rm SL} (2 , {\mathbf F}_q
[t]).
$$
Similarly we get an amalgamated structure on
$$
\Gamma_n = \biggl\{ \gamma \in {\rm SL} (2,{\mathbf F}_q [t]
[v^{-1}]) \mid \gamma \equiv \begin{pmatrix} 1 &0 \\ 0 &1
\end{pmatrix} ({\rm mod} \, n) \biggl\}.
$$
Namely $\Gamma_n = \Gamma (n) *_{\Gamma_0 (v) \cap \Gamma (n)}
\Gamma (n)$. The two injections in the amalgamated product are
the conjugate injections $\Gamma_0 (v) \cap \Gamma (n)
\hookrightarrow \Gamma (n)$ (see [Se] Ch II 1.4). Using the
Lyndon exact sequence [Se], we get the exactness of the following
sequence:
$$
H^1 (\Gamma_n , {\mathbb Z} / l) \to H^1 (\Gamma (n) , {\mathbb Z}
/ l) \oplus H^1 (\Gamma (n) , {\mathbb Z} / l) \to H^1 (\Gamma_0 (v)
\cap \Gamma (n) , {\mathbb Z} / l) \, .
$$
We consider the group $H^1 (\Gamma_n , {\mathbb Z} / l)$. An element
of this group is a homomorphism $\Gamma_n \to  {\mathbb Z} / l$.
Serre has introduced an obstruction called the congruence group
$C (G_{{\mathcal D} (S)})$ to investigate the congruence subgroup
property for a connected linear algebraic group defined over a global
field $K$ with $S$ a non-empty set of places of $K$ containing
$S_{\infty}$. Thus $C (G_{{\mathcal D} (S)})$ vanishes if and only
if the congruence subgroup property holds. Here
$$
{\mathcal D} (S) = \{ x \in K \mid v(x) \leq 1 \quad \hbox{for}
\quad v \notin S \} \, .
$$
Serre shows that for ${\rm SL} (2)$ the congruence group $C
(G_{{\mathcal D} (S)})$ is central (see [Se2] corollary of
proposition 5). Central means that $C (G_{{\mathcal D} (S)})$ is
contained in the center of the completion $\hat G$. The
completion $\hat G$ is with respect to the topology induced by
subgroups of finite index in ${\rm SL}_2 (A)$. Hence for $G =
{\rm SL} (2)$ we have $C (G_{{\mathcal D} (S)}) \simeq \pi_1
(G(S) , G(K)) := \pi_1 (G(S)) / {\rm im} (\pi_1 G(K)) \to \pi_1
(G(S))$ (see [Mo]). Moore shows that the fundamental group $\pi_1
(G(S) , G_K)$ vanishes in the special case $G = {\rm SL} (2 ,
{\mathbf F}_q (t))$ and $S = S_{\infty} \cup \{ v \}$ (see [Mo]).
As a consequence ${\rm SL} (2 , {\mathbf F}_q [t] [v^{-1}])$ has
the congruence subgroup property and the kernel of the map
$\Gamma_n \to {\mathbb Z} / l$ contains a principal congruence
subgroup. So we have ${\rm Hom} (\Gamma_n , {\mathbb Z} / l) =0$
whenever $l$ does not divide $2q$, because ${\rm PSL}_2 ({\mathbf
F}_q)$ is simple for $q > 5$. The case of ${\rm PSL}_2 ({\mathbf
F}_9)$ which is exception can be checked directly. So we get an
injection
$$
H^1 (\Gamma (n) , {\mathbb Z} / l) \oplus H^1 (\Gamma (n) ,
{\mathbb Z} / l) \to H^1 (\Gamma_0 (v) \cap \Gamma (n) ,
{\mathbb Z} / l) \eqno (*)
$$
whenever $l$ does not divide $2q$. We have an injection $H^1
(\overline{\Gamma (n)} ,  {\mathbb Z} / l) \to H^1 (\Gamma (n)$,
${\mathbb Z} / l)$. So we have shown that there exists an injection
$$
H^1 (\overline{\Gamma (n)} ,  {\mathbb Z} / l) \oplus H^1
(\overline{\Gamma (n)} ,  {\mathbb Z} / l) \to H^1
(\overline{\Gamma_0 (v) \cap \Gamma (n)} , {\mathbb Z} / l) \, .
$$
By duality the map $\alpha : \overline{\Gamma_0 (v) \cap \Gamma (n)}
\to \overline{\Gamma (n)}^{\oplus 2}$ which is induced by the
conjugate inclusions is a surjection modulo $l$ for $l$ not
dividing $2q$. $\Box$

\bigskip

\noindent {\bf Remark 2.2.2.} Note that the injectivity of $(*)$
implies that the $l$-torsion of the cokernel of $H^1 (\Gamma (n) ,
{\mathbb Z}) \oplus H^1 (\Gamma (n) ,  {\mathbb Z}) \to H^1
(\Gamma_0 (v) \cap \Gamma (n) ,  {\mathbb Z})$ vanishes.

\medskip

We have defined a ${\mathbb Q}$-valued bilinear pairing on
$\bar\Gamma$. This pairing is given by $\log \vert c_{\alpha}
(\beta) \vert$. We know that $c_{\alpha} (\beta) \in K_{\infty}$
(see [Ge-Re] p. 67). So the image of this pairing is contained in
${\mathbb Z}$. We get an injection of $\bar\Gamma \to {\rm Hom}
(\bar\Gamma , {\mathbb Z})$. We have the following commutative
diagram with bijective vertical maps:
$$
\begin{CD}
\bar\Gamma @>{\log \vert \bar c \vert}>> {\rm Hom} (\bar\Gamma ,
{\mathbb Z}) \\
@V{}VV @V{}VV  \\
\underline H_! ({\mathcal I} , {\mathbb Z})^{\Gamma} @>{}>>
{\rm Hom} (\underline H_! ({\mathcal I} , {\mathbb Z})^{\Gamma}
, {\mathbb Z}) \, .
\end{CD}
$$
The injection in the second row is defined by the following pairing
$$
(\phi_{\alpha} , \phi_{\beta})_{\mu} = \int_{E(\Gamma \backslash
{\mathcal I})} \phi_{\alpha} (e) \cdot \phi_{\beta} (e) \, n(e)^{-1} \, .
$$
Using this diagram we can control cokernel of $\bar c$.

\medskip

\noindent {\bf Lemma 2.2.3.} {\it Let $n$ be a nonzero effective
divisor on $X - \{ \infty \}$ and let $\Gamma \subset \Gamma (n)$
be a congruence subgroup. There exists a finite set of primes
$S_{\Gamma}$ which can be explicitly calculated in terms of the
graph $\Gamma \backslash {\mathcal I}$ such that for $l$ not in
$M$ the cokernel of the injection $\bar\Gamma \overset{\bar c}
{\longrightarrow} {\rm Hom} (\bar\Gamma , {\mathbb Z})$ has no
$l$-torsion.}

\medskip

\noindent {\bf Proof.} Let $T$ be a maximal tree in $\Gamma
\backslash {\mathcal I}$ and let $\{ \tilde e_1 , \ldots ,
\tilde e_g \}$ be a set of re\-presentatives for the edges
$E(\Gamma \backslash {\mathcal I}) - E(T)$ modulo orientation.
Let $v_i$ and $w_i$ denote the initial and terminal points of
$\tilde e_1$ for $i = 1 , \ldots , g$. There exists a unique
geodesic contained in $T$ which connects $w_i$ and $v_i$. Let $c_i$
be the closed path around $v_i$ obtained by this geodesic and $e_i$.
Define the following harmonic cochains $\phi_i \in \underline H_!
({\mathcal I} , {\mathbb Z})^{\Gamma}$:
$$
\phi_i (\tilde e) = \begin{pmatrix}
n(e) \, , &\hbox{if} &\tilde e \\
-n(e) \, , &\hbox{if} &\bar\tilde{\!\!e} \\
0
\end{pmatrix} \ \hbox{appears in} \ c_i \, .
$$
The inclusion $\Gamma \subset \Gamma (n)$ implies that $n(e) = 1$.
Lift $c_i$ to ${\mathcal I}$ in order to get a path from $v'_i$ to
$v''_i$ projecting to $v_i$. There exists $\alpha_i \in \Gamma$
such that $\alpha_i (v'_i) = v''_i$. In the special case where $X
- \{ \infty \}$ is the affine line, we know that $\bar\alpha_i$
form a basis for $\bar\Gamma$ (see [Re]). The fact that $n(e) = 1$
implies that $\phi_i$ forms a basis for $\underline H_!
({\mathcal I} , {\mathbb Z})^{\Gamma}$ (see [Ge-Re] proposition
3.4.5). So $\bar\alpha_i \mapsto \phi_i$ induces an isomorphism
$\bar\Gamma \to H_! ({\mathcal I} , {\mathbb Z})^{\Gamma}$. Also
we know that $\mu (e) = n(e)^{-1}$ (see [Ge-Re] 4.8). Therefore
$(\phi_i , \phi_j)_{\mu} = \sharp (c_i \cap c_j)$. Having this
formula we have an explicit matrix form $(\sharp (c_i \cap
c_j))_{i,j}$ for the map $\bar c$ with respect to the basis $\{
\bar\alpha_1 , \ldots , \bar\alpha_g \}$ for $\bar\Gamma$ and the
dual basis for ${\rm Hom} (\bar\Gamma , {\mathbb Z})$. One can
calculate the cokernel in terms of $\sharp (c_i \cap c_j)$. We
define $S_{\Gamma}$ to be the set of primes dividing $\det
(\sharp (c_i \cap c_j))_{i,j}$. We have injection modulo a prime
$l \in S_{\Gamma}$. Therefore the cokernel does not have
$l$-torsion. The graph $\Gamma \backslash {\mathcal I}$ covers
${\rm SL}_2 ({\mathbf F}_q [t]) \backslash {\mathcal I}$ with
degree $[{\rm SL}_2 ({\mathbf F}_q [t]) : \Gamma]$. So $\sharp
(c_i \cap c_j)$ is bounded by $(2 \, {\rm diam} ({\rm SL}_2
({\mathbf F}_q [t]) \backslash {\mathcal I}) + 1) [{\rm SL}_2
({\mathbf F}_q [t]) : \Gamma]$. This will help to bound primes in
$S_{\Gamma}$. $\Box$

\bigskip

\noindent {\bf Theorem 2.2.4.} {\it Let $n$ be a nonzero effective
divisor on $X - \{ \infty \}$ and let $v$ be a point on
$X - \{ \infty \}$ which does not intersect $n$. For $l$ not
dividing $2q (q_v + 1)$ and not contained in $S_{\Gamma (n)}$ we
get an injection
$$
{\rm Jac} (\bar M_{\Gamma (n)}) [l] \oplus {\rm Jac}
(\bar M_{\Gamma (n)}) [l] \to {\rm Jac} (\bar M_{\Gamma_0 (v)
\cap \Gamma (n)}) [l]
$$
which is induced by $\pi^* \circ \pi_1 + \pi'^* \circ \pi_2$.}

\medskip

\noindent {\bf Proof.} We have the following commutative diagram
with exact rows:
$$
\begin{CD}
{\scriptstyle 0} \!\! @>{}>> \!\! {\scriptstyle {\rm Hom}
(\overline{\Gamma (n)}, \mu_l)^{\oplus 2}} \!\! @>{}>> \!\!
{\scriptstyle {\rm Jac} (\bar M_{\Gamma (n)}) (C) [l]^{\oplus 2}}
\!\! @>{}>> \!\! {\scriptstyle \frac{\overline{\Gamma (n)}^{\oplus 2}}
{l \cdot \overline{\Gamma (n)}^{\oplus 2}}} \!\! @>{}>> \!\!
{\scriptstyle 0} \\
&&@V{}VV @V{}VV @V{}VV \\
{\scriptstyle 0} \!\! @>{}>> \!\! {\scriptstyle {\rm Hom}
(\overline{\Gamma_0 (v) \cap \Gamma (n)} , \mu_l)} \!\!
@>{}>> \!\! {\scriptstyle {\rm Jac} (\bar M_{\Gamma_0 (v) \cap
\Gamma (n)})(C)[l]} \!\! @>{}>> \!\! {\scriptstyle \frac{\overline{\Gamma_0 (v)
\cap \Gamma (n)}}{l \cdot \overline{\Gamma_0 (v) \cap \Gamma (n)}}}
\!\! @>{}>> \!\! {\scriptstyle 0} \, .
\end{CD}
$$
Since ${\rm Hom} (\Gamma , \mu_l) = H^1 (\Gamma , \mu_l)$ and
${\rm Hom} (\bar\Gamma , \mu_l) \to H^1 (\Gamma , \mu_l)$ is
injective, $(*)$ implies the injectivity of ${\rm Hom}
(\overline{\Gamma (n)} , \mu_l)^{\oplus 2} \to {\rm Hom}
(\overline{\Gamma_0 (v) \cap \Gamma (n)} , \mu_l)$. So it is
enough to show that $\alpha' = V \circ \pi_1 + w_v \circ V
\circ \pi_2 : \overline{\Gamma (n)}^{\oplus 2} \to \overline{\Gamma_0 (v)
\cap \Gamma (n)}$ is an injection modulo $l$ where $l$ is a
prime which does not divide $2q$. The following diagram is
commutative:
$$
\begin{CD}
\overline{\Gamma (n)}^{\oplus 2} &\hookrightarrow &{\rm Hom}
(\overline{\Gamma (n)} , {\mathbb Z})^{\oplus 2} \\
@V{V \circ \pi_1 + w_v \circ V \circ \pi_2}VV @VV{(I ,
I \circ w_v)'}V  \\
\overline{\Gamma_0 (v) \cap \Gamma (n)} &\qquad \hookrightarrow
\qquad &{\rm Hom} (\overline{\Gamma_0 (v) \cap \Gamma (n)} ,
{\mathbb Z}) \, .
\end{CD}
$$
The prime $l$ is not contained in $S_{\Gamma (n)}$. So the map
$\overline{\Gamma (n)}^{\oplus 2} \to {\rm Hom} (\overline{\Gamma
(n)} , {\mathbb Z})^{\oplus 2}$ is injective modulo $l$. We know
that ${\rm Hom} (\overline{\Gamma (n)} , {\mathbb Z})^{\oplus 2}
\to {\rm Hom} (\overline{\Gamma_0 (v) \cap \Gamma (n)} ,
{\mathbb Z})$ is injective modulo $l$. So we get injectivity of
$\overline{\Gamma (n)}^{\oplus 2} \to \overline{\Gamma_0 (v) \cap
\Gamma (n)}$. $\Box$

\bigskip

\noindent {\bf Proposition 2.2.5.} {\it Let $n$ be an effective
divisor on $X - \{ \infty \}$ and let $v$ be a point on
$X - \{ \infty \}$ which does not intersect $n$. We have a
surjection induced by $(\pi_* , \pi'_*)$:}
$$
{\rm Jac} (\bar M_{\Gamma_0 (v) \cap \Gamma (n)}) \to {\rm Jac}
(\bar M_{\Gamma (n)}) \oplus {\rm Jac} (\bar M_{\Gamma (n)})  \, .
$$

\noindent {\bf Proof.}  We will show that ${\rm Hom}
(\overline{\Gamma_0 (v) \cap \Gamma (n)} , C^*) \to {\rm Hom}
(\overline{\Gamma (n)} , C^*)^{\oplus 2}$ is surjective. We have
an injection $\overline{\Gamma (n)}^{\oplus 2} \to
\overline{\Gamma_0 (v) \cap \Gamma (n)}$. Let $L$ be a lattice in
$\overline{\Gamma_0 (v) \cap \Gamma (n)} $ such that $L \otimes
{\mathbb Q}$ be the orthogonal complement of $\overline{\Gamma
(n)}^{\oplus 2} \otimes {\mathbb Q}$ with respect to the
non-degenerate pairing on $\overline{\Gamma_0 (v) \cap \Gamma (n)}
\otimes {\mathbb Q}$. Then we have an injection from $L' = L +
\overline{\Gamma (n)}^{\oplus 2}$ into $\overline{\Gamma_0 (v)
\cap \Gamma (n)}$ with finite index. We get a surjective map ${\rm
Hom} (L' , C^*) \to {\rm Hom} (\overline{\Gamma (n)} ,
C^*)^{\oplus 2}$. This proves the surjectivity of \break ${\rm
Hom} (\overline{\Gamma_0 (v) \cap \Gamma (n)} , C^*) \to {\rm Hom}
(\overline{\Gamma (n)} , C^*)^{\oplus 2}$. $\Box$

\bigskip

Now we compute the action of  Hecke correspondences on the two copies
of $v$-old forms in order to compute the congruence module. We have a
surjection
$$
T_l ({\rm Jac} (\bar M_{\Gamma_0 (v) \cap \Gamma (n)})) \otimes
{\mathbb Q}_l \to T_l ({\rm Jac} (\bar M_{\Gamma(n)})) \otimes
{\mathbb Q}_l \oplus T_l ({\rm Jac} (\bar M_{\Gamma(n)})) \otimes
{\mathbb Q}_l \, .
$$

The cohomology group $H^1 (\bar M_{\Gamma} , {\mathbb Z}_l)$ is
the dual of the $l$-adic Tate module $T_l ({\rm Jac} (\bar
M_{\Gamma}))$. Therefore, for $l$ prime to $p$ we get an injection
$$
H^1 ({\rm Jac} (\bar M_{\Gamma (n)}) , {\mathbb Z}_l) \oplus H^1
({\rm Jac} (\bar M_{\Gamma (n)}) , {\mathbb Z}_l) \to H^1 ({\rm Jac}
(\bar M_{\Gamma_0 (v) \cap \Gamma (n)}) , {\mathbb Z}_l) \, .
$$
By theorem 2.2.4 we still have an injection modulo $l$ for $l
\nmid 2q(q_v+1)$ and not inside $S_{\Gamma (n)}$. So we have an
injection torsion-free cokernel. We shall point out that $H^1
(\bar M_{\Gamma (n)} , {\mathbb Q}_l) = H^1 ({\rm Jac} (\bar
M_{\Gamma (n)}) , {\mathbb Q}_l)$ and therefore $H^1 ({\rm Jac}$
$(\bar M_{\Gamma (n)}) , {\mathbb Q}_l)$ is a space of
automorphic forms (see theorem 1.6.1).

Let $L = H^1 ({\rm Jac} (\bar M_{\Gamma_0 (v) \cap \Gamma (n)}) ,
{\mathbb Z}_l)$  and $L' = H^1 ({\rm Jac} (\bar M_{\Gamma (n)})$,
${\mathbb Z}_l)$. Also let $V = L \otimes {\mathbb Q}_l$ and $V' =
L' \otimes {\mathbb Q}_l$ denote the associated vector spaces. The
cokernel torsion-free injection map induced between lattices in
the vector spaces $H^1 ({\rm Jac} (\bar M_{\Gamma (n)}) ,
{\mathbb Q}_l)$ and $H^1 ({\rm Jac} (\bar M_{\Gamma_0 (v) \cap
\Gamma (n)}) , {\mathbb Q}_l)$ induces an injection
$$
\beta : (V' / L')^2 \to V/L.
$$
We identify two copies of $L'$ and $V'$ with their images in $L$
and $V$ respectively. Let $A$ denote the ${\mathbb Q}_l$-vector
space generated by the image of $(L')^2$ and $B$ denote the
orthogonal subspace under the cup product. $A$ is the same as
$H^1 ({\rm Jac} (\bar M_{\Gamma_0 (v) \cap \Gamma (n)}) ,
{\mathbb Q}_l)^{v - {\rm old}}$. Since the Hecke operators $T_v$
with $v$ prime to the level are self-adjoint with respect to cup
product and cup product restricted to $A$ in non-degenerate as a
result of Poincar\'e duality on $H^1 ({\rm Jac} (\bar M_{\Gamma
(n)}) , {\mathbb Q}_l)$ we deduce that $B$ is stable under the
given Hecke operators and is the same as $H^1 ({\rm Jac} (\bar
M_{\Gamma_0 (v) \cap \Gamma (n)}) , {\mathbb Q}_l)^{v - {\rm
new}}$. We define the congruence module by the formula
$$
\Omega = ((L + A) \cap (L + B)) / L \, .
$$
One can show that $\Omega := {\rm image} (\beta) \cap \ker
(\beta') = \ker (\beta' \circ \beta)$, where the map $\beta' : V/L
\to (V'/L')^2$ is the transpose of $\beta$, induce by Poincar\'e
duality. To get congruences one should compute the congruence
module. We define the correspondence $W_v$ on $\bar M_{\Gamma_0
(v) \cap \Gamma (n)}$ by the following formula $W_v = \pi_* \pi^*
w_v - w_v$ where $\pi$ denotes the natural projection $\bar
M_{\Gamma_0 (v) \cap \Gamma (n)} \to \bar M_{\Gamma (n)} $.

\medskip

\noindent {\bf Proposition 2.2.6.} {\it The correspondence
$W_v^2 - {\rm id}$ acts on $(V' / L')^2$. We have $\Omega =
\ker (W_v^2 - {\rm id})$.}

\medskip

\noindent {\bf Proof.} We calculate the action of $W_v$ on $(V' /
L')^2$. For a cohomology class $f$ on $H^1 ({\rm Jac} (\bar
M_{\Gamma_0 (v) \cap \Gamma (n)}) , {\mathbb Q}_l)$ we have $W_v
f = \pi_* \pi^* w_v f - w_v f$. Therefore, we have $W_v w_v f =
\pi_* \pi^* f - f$. If $f$ is pull back of a cohomology class in
$H^1 ({\rm Jac} (\bar M_{\Gamma (n)}) , {\mathbb Q}_l)$ then $W_v
w_v f = q_v f$. So $W_v$ preserves the subspace $V'^2$ of $V$, and
acts as the matrix $\begin{pmatrix} T_v &q_v \\ -1 &0
\end{pmatrix}$. For $f$ in $B$ we have $\pi_* f = 0$. So $W_v$
acts on $B$ as $-w_v$. Therefore $W_v^2 - {\rm id}$ vanishes on
$B$. Now we calculate action of the correspondence $\beta' \circ
\beta$ on $V'^2$. We have
$$
\langle \pi^* f , h \rangle_{\Gamma_0 (v) \cap \Gamma (n)} =
\langle f , \pi_* h \rangle_{\Gamma (n)} \quad \hbox{for} \quad f
\in V' , \, h \in V
$$
$$
\langle w_v \pi^* f , h \rangle_{\Gamma_0 (v) \cap \Gamma (n)} =
\langle f , \pi_* w_v h \rangle_{\Gamma (n)} \quad \hbox{for}
\quad f \in V' , \, h \in V \, .
$$
The Hecke operator $T_v$ is self adjoint. Therefore for $f,g \in
V'$ we have
$$
\beta' \circ \beta (f \ g)^t = \beta' (\pi^* f + w_v \pi^* g) =
(\pi_* \pi^* f + \pi_* w_v \pi^* g \ \pi_* w_v \pi^* f + \pi_*
\pi^* g)^t \, .
$$
So the matrix of $\beta' \circ \beta$ on $V'^2$ is
$\begin{pmatrix} q_v + 1 &T_v \\ T_v &q_v + 1 \end{pmatrix}$. We
have the equality of matrices
$$
W_v^2 - {\rm id} = \begin{pmatrix} T_v &q_v \\ -1 &0
\end{pmatrix}^2 - \begin{pmatrix} 1 &0 \\ 0 &1 \end{pmatrix} =
\begin{pmatrix} -1 &T_v \\ 0 &-1 \end{pmatrix} \circ \beta' \circ
\beta \, .
$$
Since $\begin{pmatrix} -1 &T_v \\ 0 &-1 \end{pmatrix}$ acts as an
automorphism, we have $\ker (\beta' \circ \beta) = \ker (W_v^2 -
{\rm id})$.~$\Box$

\bigskip

We are interested in getting congruences between Hecke
eigen-forms. Congruence conditions should be formulated in terms
of the eigenvalues.

\medskip

\noindent {\bf Main Theorem 2.2.7.} {\it Let $n$ be a nonzero
effective divisor on $X - \{ \infty \}$ and let $v$ be a point on
$X - \{ \infty \}$ which does not intersect $n$. Let $l$ be a
prime not dividing $2q (q_v + 1)$ and not contained in $S_{\Gamma
(n)}$. Let $f$ be a Hecke eigen-form of level $\Gamma (n)$ and
$T_v f = t_v f$. If $t_v^2 \equiv (q_v + 1)^2$ mod $l$ then $f$
is congruent to a new-form of level $\Gamma_0 (v) \cap \Gamma (n)$
mod $l$.}

\medskip

\noindent {\bf Proof.}  On the two copies of old forms we have
$W_v^2 - T_v W_v + q_v = 0$. Let $r$ and $s$ denote the two roots
of $x^2 - t_v x + q_v = 0$ and let $K$ be a number field
containing $r$ and $s$ with ring of integers ${\mathcal O}_K$.
Then $f_r = f - sw_v f$ is an eigen-form of Hecke operators with
eigenvalue $r$ for $W_v$ (see [Di]). The congruence module for
$f_r$ is $(r^2 - 1)^{-1} {\mathcal O}_K / {\mathcal O}_K$. The
equality $(r^2 - 1) (s^2 - 1) = - (t_v^2 - (q_v + 1)^2)$ finishes
the argument.~$\Box$

\bigskip

\noindent {\bf Theorem 2.2.8.} {\it Let $v$ be a point on
$X - \{ \infty \}$ and let $n$ be an effective divisor on
$X - \{ \infty \}$ which does not intersect $v$. The projections
$\pi$ and $\pi'$ induce an injection of weight two and type one
modular forms
$$
M_{2,1}^2 (\Gamma (n)) (C) \oplus M_{2,1}^2 (\Gamma (n)) (C)
\to M_{2,1}^2 (\Gamma_0 (v) \cap \Gamma (n)) (C) \, .
$$}

\noindent {\bf Proof.} We have obtained a surjection
$$
{\rm Jac}
(\bar M_{\Gamma_0 (v) \cap \Gamma (n)})  \to {\rm Jac} (\bar
M_{\Gamma (n)}) \oplus {\rm Jac}(\bar M_{\Gamma (n)})
$$
induced by $(\pi_* , \pi'_*)$. So we get a surjection of tangent
spaces. This induces an injection ${\rm cot} ( {\rm Jac} (\bar
M_{\Gamma (n)})) \oplus {\rm cot} ( {\rm Jac} (\bar M_{\Gamma
(n)})) \to {\rm cot} ( {\rm Jac} (\bar M_{\Gamma_0 (v) \cap
\Gamma (n)}))$ between cotangent spaces. The cotangent space of
${\rm Jac} (\bar M_{\Gamma})$ can be cano\-nically identified
with $M_{2,1}^2 (\Gamma)(C)$. therefore, the induced map
$$
M_{2,1}^2 (\Gamma (n)) (C) \oplus M_{2,1}^2 (\Gamma (n))(C) \to
M_{2,1}^2 (\Gamma_0 (v) \cap \Gamma (n)) (C)
$$
is an injection.$\Box$

\bigskip

We know that the vector spaces $M_{2,1}^2 (\Gamma)(C)$ have
${\mathbf F}_q$ structure. We are curious if the above injection
is also defined over ${\mathbf F}_q$. To show this we have to use
the language of harmonic cochains.

\medskip

\noindent {\bf Proposition 2.2.9.} {\it $M_{2,1}^2 (\Gamma (n))
(C) \oplus M_{2,1}^2 (\Gamma (n)) (C) \to M_{2,1}^2 (\Gamma_0 (v)
\cap \Gamma (n)) (C)$ can be reduced to ${\mathbf F}_q$ to induce
an injection}
$$
M_{2,1}^2 (\Gamma (n)) ({\mathbf F}_q) \oplus M_{2,1}^2 (\Gamma
(n)) ({\mathbf F}_q) \to M_{2,1}^2 (\Gamma_0 (v) \cap \Gamma (n))
({\mathbf F}_q) \, .
$$

\noindent {\bf Proof.} The equalities $M_{2,1}^2 (\Gamma) (C) =
\underline H_! ({\mathcal I} , C)^{\Gamma}$ and $M_{2,1}^2
(\Gamma) ({\mathbf F}_q) = \underline H_! ({\mathcal I} ,
{\mathbf F}_q)^{\Gamma}$ induce a map $\underline H_! ({\mathcal
I} , C)^{\Gamma (n)} \oplus \underline H_! ({\mathcal I} ,
C)^{\Gamma (n)} \to \underline H_! ({\mathcal I} , C)^{\Gamma_0
(v) \cap \Gamma (n)}$. The first copy maps by the natural
injection $\underline H_! ({\mathcal I} , C)^{\Gamma (n)} \to
\underline H_! ({\mathcal I} , C)^{\Gamma_0 (v) \cap \Gamma (n)}$
and the second map is given by the first map composed with the
involution $w_v$ of $ \underline H_! ({\mathcal I} , C)^{\Gamma_0
(v) \cap \Gamma (n)}$. We can restrict to ${\mathbf F}_q$ and
obtain an injective map $\underline H_! ({\mathcal I} , {\mathbf
F}_q)^{\Gamma (n)} \oplus \underline H_! ({\mathcal I} , {\mathbf
F}_q)^{\Gamma (n)} \to \underline H_! ({\mathcal I} , {\mathbf
F}_q)^{\Gamma_0 (v) \cap \Gamma (n)}$. This means that by
restriction we get $M_{2,1}^2 (\Gamma (n)) ({\mathbf
F}_q)^{\oplus 2} \to M_{2,1}^2 (\Gamma_0 (v) \cap \Gamma (n))
({\mathbf F}_q)$. $\Box$

\bigskip

After tensoring $M_{2,1}^2 (\Gamma (n)) ({\mathbf F}_q)^{\oplus
2} \to M_{2,1}^2 (\Gamma_0 (v) \cap \Gamma (n)) ({\mathbf F}_q)$
with ${\mathbf F}_q [t]$ or with the ring of integers $A$ of any
function field, we get a cokernel torsion-free injection. In
particular we have a cokernel torsion-free injection of ${\mathbf
F}_q [t]$-modules
$$
M_{2,1}^2 (\Gamma (n)) ({\mathbf F}_q [t])^{\oplus 2} \to
M_{2,1}^2 (\Gamma_0 (v) \cap \Gamma (n)) ({\mathbf F}_q [t]) \, .
$$
Unfortunately we don't have a ${\mathbf F}_q [t]$-valued pairing
on these spaces. We can try to construct a pairing on $M_{2,1}^2
(\Gamma) ({\mathbf F}_q)$ using the pairing on $\underline H_!
({\mathcal I} , {\mathbf F}_q)^{\Gamma}$. But it is not sensible
to search for congruences between ${\mathbf F}_q [t]$-valued
Drinfeld modular forms using such a pairing.

\subsection{Towers of congruences}\label{ssec2.3}

In order to get congruences between Hecke algebras we need to
construct congruences between level $\Gamma_0 (rv) \cap \Gamma
(n)$ and level $\Gamma_0 ((r + 2) v) \cap \Gamma (n)$. The first
thing we need is an injection result on the $l$-torsion of
Jacobians. The main reference in [Wi].

Conjugation by $\tau_{rv} = \begin{pmatrix} \pi_v^r &0 \\ 0 &1
\end{pmatrix}$ gives an isomorphism between the abelianizations
$\overline{\Gamma_0 (rv) \cap \Gamma (n)} \to \overline{\Gamma^0
(rv) \cap \Gamma (n)}$. Conjugation of $\Gamma_0 (r,v)$ by
$\begin{pmatrix} 0 &1 \\ -1 &0
\end{pmatrix}$ gives another isomorphism. Composing one with the
inverse of the other induces an involution $w_v :
\overline{\Gamma_0 (rv) \cap \Gamma (n)} \to \overline{\Gamma_0
(rv) \cap \Gamma (n)}$. This induces an involution $w'_{rv}$ of
${\rm Hom} (\overline{\Gamma_0 (rv) \cap \Gamma (n)} , C^*)$ and
an involution $w''_{rv}$ of ${\rm Jac} (\bar M_{\Gamma_0 (rv)
\cap \Gamma (n)})$. Now we have three maps from $\bar M_{\Gamma_0
((r+2)v) \cap \Gamma (n)} \to \bar M_{\Gamma_0 (rv) \cap \Gamma
(n)}$. The following noncommutative diagram helps to visualize
these maps
$$
\begin{CD}
\bar M_{\Gamma_0 ((r+2)v) \cap \Gamma (n)} @>{\pi_{00}}>>
\bar M_{\Gamma_0 ((r+1)v) \cap \Gamma (n)} @>{\pi_0}>>
\bar M_{\Gamma_0 (rv) \cap \Gamma (n)} \\
@V{w_{(r+2)v}}VV @A{w_{(r+1)v}}AA @V{w_{rv}}VV \\
\bar M_{\Gamma_0 ((r+2)v) \cap \Gamma (n)} @>{\pi_{00}}>>
\bar M_{\Gamma_0 ((r+1)v) \cap \Gamma (n)} @>{\pi_0}>>
\bar M_{\Gamma_0 (rv) \cap \Gamma (n)} \, .
\end{CD}
$$

Here $\pi_0$, $\pi_{00}$ denote the natural projections. The
three maps are $\alpha = \pi_0 \pi_{00}$ and $\beta = w_{rv}
\pi_0 w_{(r+1)v} \pi_{00}$ and $\gamma = w_{rv} \pi_0 \pi_{00}
w_{(r+2)v}$. Let $\pi_i$ denote projection to the $i$-th
component.

\medskip

\noindent {\bf Theorem 2.3.1.} {\it Let $n$ be an effective
divisor on $X - \{ \infty \}$ and let $v$ be a point on  $X - \{
\infty \}$ which does not intersect $n$. For $l$ not dividing
$q(q_v + 1)$ which is not contained in $S_{\Gamma_0 (rv) \cap
\Gamma (n)}$ we get an exact sequence of $l$-torsion of Jacobians:
$$
0 \longrightarrow {\rm Jac} (\bar M_{\Gamma_0 (rv) \cap \Gamma
(n)}) [l] \overset{\zeta_1}{\longrightarrow} {\rm Jac} (\bar
M_{\Gamma_0 ((r+1)v) \cap \Gamma (n)}) [l]^{\oplus 2} $$
$$
\overset{\zeta_2}{\longrightarrow}  {\rm Jac} (\bar M_{\Gamma_0
((r+2)v) \cap \Gamma (n)}) [l]
$$
where $\zeta_1 = (\pi_{00}^* , - w_{(r+1)v}^* \pi_0^* w_{rv}^*)$
and $\zeta_2 = \pi_{00}^* + w_{(r+2)v}^* \pi_{00}^*
w_{(r+1)v}^*$.}

\medskip

\noindent {\bf Proof.} See [Wi] lemma 2.5 for more details. Let
$B_0$, $B^0$ and $\Delta_{rv}$ be given by
$$
B_0 = \Gamma_0 ((r+1)v) \cap \Gamma (n) / \Gamma_0 ((r+1)v) \cap
\Gamma (v) \cap \Gamma (n) \, ,
$$
$$
B^0 = \Gamma_0 (rv) \cap \Gamma^0 (v) \cap \Gamma (n) / \Gamma_0
(rv) \cap \Gamma (v) \cap \Gamma (n) \, ,
$$
$$
\Delta_{rv} = \Gamma_0 (rv) \cap \Gamma (n) / \Gamma_0 ((r+1)v)
\cap \Gamma (v) \cap \Gamma (n) \, .
$$
Then $\Delta_{rv} \simeq {\rm SL}_2 ({\mathbf F}_{q_v})$ for
$r=0$ and is of order a power of $q_v$ for $r > 0$. The groups
$B_0$ and $B^0$ generate $\Delta_{rv}$. The vanishing of $H^2
({\rm SL}_2 ({\mathbf F}_{q_v}) , {\mathbb Z} / l)$ can be
checked by restriction to the Sylow $l$-subgroup which is cyclic.
The following isomorphisms
$$
\lambda_0 : H^1 (\Gamma_0 ((r+1)v) \cap \Gamma (n) , {\mathbb Z}
/ l) \overset{\simeq}{\longrightarrow} H^1 (\Gamma_0 ((r+1)v)
\cap \Gamma (v) \cap \Gamma (n) , {\mathbb Z} / l)^{B_0} \, ,
$$
$$
\lambda^0 : H^1 (\Gamma_0 ((r+1)v) \cap \Gamma (n) , {\mathbb Z}
/ l) \overset{\simeq}{\longrightarrow} H^1 (\Gamma_0 ((r+1)v)
\cap \Gamma (v) \cap \Gamma (n) , {\mathbb Z} / l)^{B^0} \, ,
$$
$$
H^1 (\Gamma_0 (rv) \cap \Gamma (n) , {\mathbb Z} / l)
\overset{\simeq}{\longrightarrow} H^1 (\Gamma_0 ((r+1)v) \cap
\Gamma (v) \cap \Gamma (n) , {\mathbb Z} / l)^{\Delta_{rv}}
$$
is induced by Inflation-Restriction exact sequences. We have the
inclusion $H^1 (\Gamma)^{B^0} \cap H^1 (\Gamma)^{B_0} \subset H^1
(\Gamma)^{\Delta_{rv}}$. This implies exactness the following
sequence:
$$
0 \to H^1 (\Gamma_0 (rv) \cap \Gamma (n), {\mathbb Z} / l)
$$
$$
\overset{\zeta_1}{\longrightarrow} H^1 (\Gamma_0 ((r+1)v) \cap
\Gamma (n) , {\mathbb Z} / l) \oplus H^1 (\Gamma_0 (rv) \cap
\Gamma^0 (v) \cap \Gamma (n) , {\mathbb Z} / l)
$$
$$
\overset{\zeta_2}{\longrightarrow} H^1 (\Gamma_0 ((r+1)v) \cap
\Gamma (v) \cap \Gamma (n) , {\mathbb Z} / l)
$$
where $m_1 = ({\rm res}_0 , { \rm res}^0)$ and $m_2 = \lambda_0 +
\lambda^0$. By appropriate conjugation one gets the following
exact sequence:
$$
0 \to H^1 (\Gamma_0 (rv) \cap \Gamma (n) , {\mathbb Z} / l)
$$
$$
\overset{\zeta_1}{\longrightarrow} H^1 (\Gamma_0 ((r+1)v) \cap
\Gamma (n) , {\mathbb Z} / l)^{\oplus 2}
\overset{\zeta_2}{\longrightarrow} H^1 (\Gamma_0 ((r+2)v) \cap
\Gamma (n) , {\mathbb Z} / l) \, .
$$
As in Theorem 2.2.1 this gives a dual sequence which is exact
modulo $l$ \
$$
\Gamma_0 ((r+2)v) \cap \Gamma (n) \to \Gamma_0
((r+1)v) \cap \Gamma (n)^{\oplus 2} \to \Gamma_0 (rv) \cap \Gamma
(n) \to 0 \, .
$$
As in the proof of Theorem 2.2.4 we can use the Gekeler-Reversat
uniformization to get information about $l$-torsion of Jacobians.
The condition that $l$ is not in $S_{\Gamma_0 (rv) \cap \Gamma
(n)}$ is crucial to get injection of torsion of Jacobians. The
same argument gives us the exact sequence we are seeking for.
$\Box$

\bigskip

\noindent {\bf Theorem 2.3.2.} {\it Let $n$ be an effective
divisor on $X - \{ \infty \}$ and let $v$ be a point on  $X - \{
\infty \}$ which does not intersect $n$. For $l$ not dividing
$q(q_v + 1)$ which is not contained in $S_{\Gamma_0 (rv) \cap
\Gamma (n)}$, the map $\alpha^* \pi_1 + \beta^* \pi_2 + \gamma^*
\pi_3$ induces an injection}
$$
{\rm Jac} (\bar M_{\Gamma_0 (rv) \cap \Gamma (n)}) [l]^{\oplus 3}
\to {\rm Jac} (\bar M_{\Gamma_0 ((r+2)v) \cap \Gamma_0 (n)}) [l]
\, .
$$

\medskip

\noindent {\bf Proof.} Let $\pi'_0 = w_{rv}^* \pi_0 w_{(r+1)v}$
and $\pi'_{00} = w_{(r+1)v} \pi_{00} w_{(r+2)v}$, then we have
$\pi'_0 \cdot \pi_{00} = \pi_0 \cdot \pi'_{00}$. Therefore two of
the four maps $\pi_0 \cdot \pi_{00}$, $\pi'_0 \cdot \pi_{00}$,
$\pi_0 \cdot  \pi'_{00}$ and $\pi'_0 \cdot \pi'_{00}$ below
coincide
$$
{\rm Jac} (\bar M_{\Gamma_0 (rv) \cap \Gamma (n)})[l] \
\begin{matrix} \overset{\pi'^*_0}{\rightarrow} \\
\overset{\pi^*_0}{\rightarrow} \end{matrix} \ {\rm Jac} (\bar
M_{\Gamma_0 ((r+1)v) \cap \Gamma_0 (n)})[l] \ \begin{matrix}
\overset{\pi'^*_{00}}{\rightarrow} \\
\overset{\pi^*_{00}}{\rightarrow} \end{matrix} \ {\rm Jac} (\bar
M_{\Gamma_0 ((r+2)v) \cap \Gamma_0 (n)})[l]
$$
and we get injection of the three copies which are left. $\Box$

\bigskip

Let $L = H^1 (\bar M_{\Gamma_0 (rv) \cap \Gamma (n)} , {\mathbb
Z}_l)$ and $L' = H^1 (\bar M_{\Gamma_0 ((r+2)v) \cap \Gamma (n)},
{\mathbb Z}_l)$. Also let $V = L \otimes {\mathbb Q}_l$ and $V' =
L' \otimes {\mathbb Q}_l$ denote the associated vector spaces. We
have obtained an injection $\mu : (V'/L')^3 \to V/L$. We identify
three copies of $L'$ and $V'$ with their images in $L$ and $V$
respectively. Let $A$ denote the ${\mathbb Q}_l$-vector space
generated by the image of $(L')^3$ and $B$ denote the orthogonal
subspace under the cup product. We define the congruence module
by the formula
$$
\Omega = ((L+A) \cap (L+B)) / L \, .
$$
One can show that $\Omega := {\rm image} (\mu) \cap \ker (\mu') =
\ker (\mu' \circ \mu)$, where the map $\mu' : V/L \to (V'/L')^3$
is the transpose of $\beta$, induce by Poincar\'e duality. To get
congruences one should compute the congruence module. First we
calculate action of the correspondence $\mu' \circ \mu$ on
$V'^3$. For $f \in V'$ and $h \in V$ we have
$$
\langle \pi_{00}^* \pi_0^* f , h \rangle_{\Gamma_0 ((r+2)v) \cap
\Gamma (n)} = \langle f , \pi_{0*} \pi_{00*} h \rangle_{\Gamma_0
(rv) \cap \Gamma (n)},
$$
$$
\langle w_{(r+2)v} \pi_{00}^* w_{(r+1)v} \pi_0^* f , h
\rangle_{\Gamma_0 ((r+2)v) \cap \Gamma (n)} =\hspace{2 in}
$$
$$
\hspace{2 in} \langle f , \pi_{0*} w_{(r+1)v} \pi_{00*} w_{(r+2)v}
h \rangle_{\Gamma_0 (rv) \cap \Gamma (n)},
$$
$$
\langle w_{(r+2)v} \pi_{00}^* \pi_0^* w_{rv} f , h
\rangle_{\Gamma_0 ((r+2)v) \cap \Gamma (n)} = \langle f , w_{rv}
\pi_{0*} \pi_{00*} w_{(r+2)v} h \rangle_{\Gamma_0 (rv) \cap
\Gamma (n)} \, .
$$
So the matrix of the action $\mu' \circ \mu$ on $V'^3$ is given by
$$
\begin{pmatrix}
q_v (q_v + 1) &w_{rv} \, T_v \cdot q_v &w_{rv} \, T_{2v} \\
T_v \, w_{rv} \cdot q_v &q_v (q_v + 1) &w_{rv} \, T_v \cdot q_v \\
T_{2v} \, w_{rv} &T_2 \, vw_{rv} \cdot q_v &q_v (q_v + 1)
\end{pmatrix} .
$$
The Hecke operators $T_v$ and $T_{2v}$ are self adjoint and defined by
$$
T_v = w_{rv} \, \pi_{0*} \, w_{(r+1)v} \, \pi_0^* \, w_{rv}
$$
$$
T_{2v} = \pi_{0*} \, \pi_{00*} \, w_{(r+2)v} \, \pi_{00}^* \,
\pi_0^* \, w_{rv} \, .
$$
This is enough material to get congruences between Hecke
algebras. We do not intend to enter the realm of Hecke algebras
in this manuscript.

\section*{Acknowledgements}

For this research I have benefited conversations with A. Borel,
A. Rajaei, P. Sarnak. I specially thank J. K. Yu who criticized
the manuscript in detail. His questions made my mind much more
clear. I would also like to thank A. Genestier who was so kind to
give some comments on background material. He also pointed out
that he could get the corresponding Ihara-type result directly
using Hecke algebras. This paper is based on my thesis presented
years ago to Princeton University. To my advisor A. Wiles goes my
most sincere thanks.

\bibliography{Bibliography}

Sharif University of Technology, e-mail: rastegar@sharif.edu
\\Institut des Hautes Etudes Scientifiques, e-mail:
rastegar@ihes.fr

\end{document}